\documentstyle[12pt]{article}

\newcommand{\sect}[1]{\setcounter{equation}{0}\section{#1}}
\newcommand{\subsect}[1]{\subsection{#1}}

\font\black=msbm10 scaled\magstep1
\def\field #1{\hbox{{\black #1}}}
\def\Re{{\hbox{{\field R}}}}
\def\be{\begin{equation}}
\def\ee{\end{equation}}
\def\bea{\begin{eqnarray}}
\def\eea{\end{eqnarray}}

\def\om{\Omega}
\def\k{\omega}

\def\cext{\alpha}

\def\cextiif{\alpha^{F}}
\def\cextiil{\alpha^{L}}
\def\cextiii{\beta}

\def\G{{\cal G}}

\def\propi{3.1}
\def\propii{3.2}
\def\teor{4.1}

\begin{document}
\baselineskip=18pt

\thispagestyle{empty}
\rightline{DAMTP 96-86, revised October 1997}

\vspace{2cm}
\begin{center}
\begin{large}
{\bf  Central Extensions of the Quasi-orthogonal Lie algebras}
\end{large}
\\[0.4cm]

J. A. de Azc\'arraga$^\dagger$
\footnote[4]{St. John's College Overseas Visiting Scholar.}
\footnotemark[5], F. J. Herranz$^\star$,
J. C. P\'erez Bueno$^\dagger$
\footnote[5]
{On sabbatical (J.A.) leave and on leave of absence (J.C.P.B.)
from Departamento de F\'{\i}sica Te\'orica and IFIC
(Centro Mixto Univ. de Valencia-CSIC) E--46100 Burjassot (Valencia), Spain.}
and M. Santander$^\ddagger$
\\ [0.3cm]
{\it $^\dagger$ Department of Applied Mathematics and Theoretical Physics,
\\Silver St., Cambridge, CB3-9EW, UK}
\\[0.3cm]
{\it
$^\star$
Departamento de F\'{\i}sica, E.U. Polit\'ecnica
\\
Universidad de Burgos,
E--09006 Burgos, Spain
\\[0.3cm]
$^{\ddagger}$ Departamento de F\'{\i}sica Te\'orica,
Universidad de Valladolid
\\
E--47011, Valladolid, Spain
}
\end{center}

\begin{abstract}
We determine the central extensions of a whole family of Lie algebras,
obtained
by the method  of graded contractions from
$so(N+1)$, $N$ arbitrary. All the inhomogeneous orthogonal and
pseudo-orthogonal algebras are members of this family, as well as a large
number of other non-semisimple algebras, all of which have at least a
semidirect structure (in some cases two or more). The dimensions of their
second cohomology groups
$H^2(\G,\Re)$ and the explicit expression of their central extensions are
given.
\end{abstract}


\sect {Introduction}

This paper is devoted to investigate the Lie algebra cohomology groups of a
large class of algebras, the so-called orthogonal Cayley-Klein (CK)
family of algebras. These algebras may be obtained by a sequence of
ordinary contractions starting from $so(N+1)$ or from $so(p,q)$,
and can be described by using the
alternative  method of graded contractions
\cite{MonPat,MooPat}.

The problem of finding the cohomology groups for this
family of algebras is primarily of mathematical interest, but it is not devoid
of a physical one since the CK algebras include all kinematical algebras of
physical relevance.
There are (at least) three main areas where central
extensions play a r\^ole in physics. First, the existence of a
non-trivial cohomology group is  associated with projective representations, a
fact discussed first in general by Bargmann
\cite{Bar} but which has its roots in the work of Weyl
\cite{WEYL} and in the classic paper of Wigner \cite{WIGNER}.
Second, in the Kirillov-Kostant-Souriau theory, homogeneous symplectic
manifolds
under a group appear as the orbits of the coadjoint representation of either
the group itself or of a central extension. Third, if a group is considered as
the invariance group of a given physical theory, the most general Lagrangian
which leads to invariant equations of motion is not necessarily
a strictly invariant Lagrangian, but a quasi-invariant one;
again this is linked to the central extensions of the group
\cite{Levi}.

The second cohomology group is
trivial for semisimple Lie algebras. For
some specific non-semisimple
algebras algebras it has also been studied ({\it e.g.}, for Euclidean,
Poincar\'e and
Galilei algebras in $N$ dimensions). It is well known \cite{Bar,LL}
that the Galilei group admits projective  representations, and the Lie algebra
statement $H^2(\G(3+1),\Re)=\Re$ has its counterpart in the fact that the group
can be centrally extended by $U(1)$ (the phase group involved in the ray
\cite{Bar} representations).
The cohomology of these algebras in low dimensions is also known.
However, some recent rediscoveries (as for the (2+1) Galilei algebra)
suggest that part of these results have become buried in the
literature.

Moving from the semisimple $so(p,q)$ algebras
by successive
contractions it is found that, as a rule,
the cohomology groups of the resulting algebras become larger.
However, inhomogeneity is not {\it per se} a sufficient condition for the
existence of non-trivial cohomology, as exhibited by the contrast between the
ten-dimensional Galilei and Poincar\'e algebras.
Moreover, as far as we
know, there is no systematic study of the second cohomology groups covering a
large family of algebras, so that the present study may help to describe
the relationship between cohomology and contractions.
It should be noted, however, that the dimension of the Lie algebra
cohomology groups does not need to be the same as the dimension of the
$U(1)$-valued group cohomology.
This is because when moving from the Lie algebras to Lie groups the topology
of the groups comes into play and this may reduce the dimension of
the different
cohomology groups (for a general statement see \cite{Bar} and
\cite{VanEst}).
In  Bargmann's terminology \cite{Bar} (see \cite{TUYNMAN} for an outlook),
it is not always possible to extend a local exponent to
the whole group (it is possible
if the group is simply connected, hence the r\^ole of the universal covering
group).

This paper is organised as follows.
The algebras in the orthogonal CK family are introduced in sec. 2,
as well as the way
the relevant kinematical Lie algebras are included in it.
Sec. 3 gives a reasonably complete account of the procedure used to find the
central extensions of the orthogonal CK algebras (including the
characterisation of trivial and equivalent extensions).
The paper allows to skip this section completely if the reader is not
interested in the computational details.
In sec. 4 we state the final solution and provide a method to compute the
dimension of the second cohomology group for any CK algebra in any dimension.
The method is illustrated in sec. 5, where all central extensions for
the lower dimensional CK orthogonal algebras are explicitly described,
and some comments on the {\it group} cohomology
for the CK quasi-orthogonal groups are made.  Some concluding
remarks and prospects for future work close the paper.


\sect {The CK family of quasi-orthogonal algebras \break
as graded contractions of $so(N+1)$}

Let us first introduce the set of algebras we are going to deal with
in connection with graded contraction theory \cite{MonPat,MooPat},
which is a convenient tool for our aims in this paper.
The starting point is the real Lie algebra
$so(N+1)$ with
$N(N+1)/2$ generators $\om_{ab}$ $(a,b=0,1,\dots, N, \  a<b)$. Its
non-zero Lie brackets are
\be
[\om_{ab}, \om_{ac}] =  \om_{bc}\quad, \quad
[\om_{ab}, \om_{bc}] = -\om_{ac}\quad, \quad
[\om_{ac}, \om_{bc}] =  \om_{ab}\quad, \quad  a<b<c
\label{ca}
\ee
(commutators involving four different indices are zero). This Lie
algebra can be endowed with a fine  grading group ${\bf
{Z}}_2^{\otimes N}$.
The corresponding graded contractions
of $so(N+1)$ constitute a large set of Lie algebras, depending on $2^N$--$1$
real contraction parameters \cite{HerSanGSGC}, which include from  the
simple Lie algebras $so(p,q)$ (when all parameters are different from
zero) to the Abelian algebra at the opposite case, when all parameters are
equal to zero.
Properties associated with the simplicity of the algebra are lost
at some point beyond the simple algebras in the contraction lattice, yet there
exists a particular subset or family of these graded contractions,
nearer to the simple ones, which essentially preserve these
properties, and may therefore be called \cite{Ros} `quasi-simple' algebras.
This family \cite{HMOS} encompasses the
pseudo-orthogonal algebras (the $B_l$ and $D_l$ Cartan series)
as well as their nearest non-simple contractions; collectively,
all these algebras are called {\it quasi-orthogonal}. In this paper we will
deal exclusively with this family\footnote{The study of the cohomology of the
general
graded contractions of $so(N+1)$ could be performed  similarly;
we shall nevertheless restrict ourselves to the CK family.
This means, for instance, that the `completely contracted' (and hence Abelian)
algebra, for which dim$H^2(\G,\Re)=\left({\hbox{dim}\,\G\atop 2}\right)$, is
not included.}
to be defined precisely below. Its members are called Cayley--Klein algebras
since they
are exactly the family of motion algebras of the geometries of a real space
with a projective metric in the Cayley--Klein sense \cite{Sommer,Yas}.
If the graded contraction
procedure is implicitly assumed, the algebras in the CK family may be
referred to as contractions of the compact
$so(N+1)$, although in this sentence the word `contraction' must be
adequately understood and the compactness is lost.

The set of CK Lie algebras depends on $N$ real coefficients
$\k_1,\dots,\k_N$ which codify in a convenient way the information on the
Lie algebra structure. In terms of
\be
\k_{ab}:=\k_{a+1}\k_{a+2}\cdots\k_b\quad\quad (a,b=0,1,\dots,N
\quad,\quad a<b ) \quad , \quad \k_{aa}:=1\label{cc}
\ee
for which we remark the relations:
\bea
&&\k_{ac}=\k_{ab}\k_{bc}\quad,\quad a \leq b \leq c  \qquad.
\k_{a}=\k_{a-1\, a}\quad,\quad a=1,\dots,N\quad.
\label{ce}
\eea
the independent  non-vanishing commutation relations in the
contracted CK Lie algebra are given (cf. (\ref{ca})) by
\be
[\om_{ab}, \om_{ac}] =  \k_{ab}\om_{bc}\quad, \quad
[\om_{ab}, \om_{bc}] = -\om_{ac}\quad, \quad
[\om_{ac}, \om_{bc}] =  \k_{bc}\om_{ab}\quad, \quad  a<b<c \quad.
\label{cb}
\ee
This CK or quasi-orthogonal algebra will  be denoted as
$so_{\k_1,\dots,\k_N}(N+1)$, making explicit the parameters $\k_i$.
This has a (vector) representation by
$(N+1)\times (N+1)$ real matrices, given by
\be
\om_{ab}=-\k_{ab}e_{ab}+e_{ba}
\label{ccd}
\ee
where $e_{ab}$ is the matrix with a single  non-zero entry, 1,
in the row $a$ and column $b$.
For $\k_1=\ldots=\k_N=1$, we recover the compact $so(N+1)$
algebra.
Since each coefficient $\k_a$ may take positive,
negative or zero values, and by means of a simple rescaling of the
initial generators it can be reduced to the standard values of $1$, $-1$
or $0$, it is clear that the family
$so_{\k_1,\dots,\k_N}(N+1)$   includes $3^N$ Lie algebras. Some of
these can be isomorphic;   for instance,
\be
so_{\k_1,\k_2,\dots,\k_{N-1},\k_N}(N+1)\simeq
so_{\k_N,\k_{N-1},\dots,\k_2,\k_1}(N+1).
\label{cf}
\ee

The family $so_{\k_1,\dots,\k_N}(N+1)$ of CK algebras
includes algebras of physical interest
\cite{Tesis}. The structure of these algebras can be
characterised by two statements:

\noindent
{\bf $\bullet$} When {\it all} constants $\k_a\ne 0$,   the algebra
$so_{\k_1,\dots,\k_N}(N+1)$ is a simple (barring the special $N=3$ case)
real Lie algebra in the Cartan series $B_l$ or $D_l$, and it is
isomorphic to a pseudo-orthogonal algebra $so(p,q)$ with $p+q=N+1$.

\noindent
{\bf $\bullet$} If a constant $\k_a = 0$, for $a=1, \dots, N$, the
resulting algebra $so_{\k_1,\dots,\k_a=0,\dots,\k_N}(N+1)$,    has the
semidirect structure (see eq. (\ref{cb}))
\be
so_{\k_1,\dots,\k_{a-1},\k_a=0,\k_{a+1},\dots,\k_N}(N+1)
\equiv t \odot
( so_{\k_1,\dots,\k_{a-1}}(a) \oplus
so_{\k_{a+1},\dots,\k_N}(N+1-a))\quad,
\label{cg}
\ee
where $t$ is an Abelian subalgebra
$\hbox{dim}\,t=a(N+1-a)$, and  the remaining subalgebra is a direct
sum.  The three subalgebras appearing in (\ref{cg}) are generated by:
\bea
&&t=\langle \om_{ij};\ i=0,1,\dots,a-1;\ j=a,a+1,\dots,N\rangle\quad,\cr
&&so_{\k_1,\dots,\k_{a-1}}(a)
=\langle \om_{ij};\ i,j=0,1,\dots,a-1\rangle\quad,\cr
&&so_{\k_{a+1},\dots,\k_N}(N+1-a)
=\langle \om_{ij};\ i,j=a,a+1,\dots,N\rangle \quad.
\label{ccg}
\eea
If either $a=1$ or
$a=N$, then the direct sum in  (\ref{cg}) has a single summand. This
decomposition as a semidirect  sum is true irrespective of whether the
remaining constants $\k_i$ are equal to zero or not.

The structure behind this decomposition can be described visually by setting
the generators in a triangular array.
The generators spanning the subspace $t$ are those inside the rectangle, while
the subalgebras $so_{\k_1,\dots,\k_{a-1}}(a)$ and
$so_{\k_{a+1},\dots,\k_N}(N+1-a)$ correspond to
the two triangles to the left and below
the rectangle respectively,

\bigskip
\noindent\hskip 1.8truecm
\begin{tabular}{cccc|cccc}
$\om_{01} $&$ \om_{02} $&$\ldots$&
$\om_{0\, a-1} $&
 $\om_{0a}$&$\om_{0\, a+1}$&
$\ldots$&$\om_{0N}$\\
 &$ \om_{12} $&$\ldots$&$\om_{1\, a-1} $& $\om_{1a}$&$\om_{1\, a+1}$&
$\ldots$&$\om_{1N}$\\
 &&$\ddots $&$\vdots$&  $\vdots$&$\vdots$&
$ $&$\vdots$\\
 & &$ $&$\om_{a-2\,a-1}$&  $\om_{a-2\,a}$&$\om_{a-2\,a+1}$&
$\ldots$&$\om_{a-2\,N}$\\
 & &$ $& &  $\om_{a-1\,a}$&$\om_{a-1\,a+1}$&
$\ldots$&$\om_{a-1\,N}$\\
\cline{5-8}
 & &$ $&\multicolumn{1}{c}{\,}&    $ $&$\om_{a\,a+1}$&
$\ldots$&$\om_{a  N}$\\
 & &$ $&\multicolumn{1}{c}{\,}& $ $ &$ $&
$\ddots$&$\vdots$\\
 & &$ $&\multicolumn{1}{c}{\,}& $ $&$ $&
$ $&$\om_{N-1\,N}$\\
\end{tabular}

\bigskip
\noindent
The Abelian rectangle is reduced to the row $\om_{0i}\,,\,i=1,...,N$ for
$\k_1=0$ or to the column $\om_{iN}\,,\,i=0,1,...,N-1$ when $\k_N=0$.

Let us consider some specially
relevant algebras in the CK family.
All $\k$'s are assumed to be different from zero,
unless otherwise explicitly stated.

\noindent
{\bf (1)}
$\k_a\ne 0$ $\forall a$. Here $so_{\k_1,\dots,\k_N}(N+1)$ is a
pseudo-orthogonal algebra $so(p,q)$ with $p+q=N+1$. The
matrix representation (\ref{ccd}) of this algebra generates a group of
linear transformations on $N+1$ real variables which
leave invariant a quadratic form $g$, with matrix
$\hbox{diag} (1,\k_{01},\k_{02},\dots,\k_{0N})$. The signature of this
quadratic form is
the number of positive and negative terms in the sequence
$(1,\k_{1},\k_{1}\k_{2},\dots,\k_{1}\k_{2}\dots\k_{N})$, so that each $\k_i$
governs the relative sign of two consecutive diagonal elements in the
metric matrix, the element where $\k_i$ appears for the first time and
the previous one.

\noindent
{\bf (2)} $\k_1=0$.  The $so_{0,\k_2,\dots,\k_N}(N+1)$ algebras
are the usual pseudo-orthogonal inhomogeneous ones,
with a semidirect sum structure given by
$$
so_{0,\k_2,\dots,\k_N}(N+1)
\equiv t_N\odot  so_{\k_2,\dots,\k_N}(N)\equiv iso(p,q)  \quad,\quad
p+q=N\quad,
$$
where $so_{\k_2,\dots,\k_N}(N)$ acts on
$t_N$ through the vector representation.
The Euclidean algebra $iso(N)$ appears once with
$(\k_1,\k_2,\dots,\k_N)={(0,1,\dots,1)}$. The
Poincar\'e algebra  $iso(N-1,1)$ is reproduced several times, {\it e.g.} for
$(\k_1,\k_2,\dots,\k_N)=(0,-1,1,\dots,1)$.

\noindent
{\bf (3)} $\k_1=\k_2=0$.  These algebras have two different semidirect
sum structures (cf. (\ref{cg})).
The one associated to the vanishing of $\k_1$ is
$$
t_N \odot
so_{\k_2=0,\k_{3}, \dots,\k_N }(N)
$$
while the structure associated to the vanishing of $\k_2$ is
$$t_{2N-2} \odot  \left( so_{\k_1=0}(2) \oplus
so_{\k_{3},\dots,\k_N}(N-1) \right)
$$
The first structure can also be seen as a `twice-inhomogeneous'
pseudo-orthogonal algebra
\be
so_{0,0,\k_3,\dots,\k_N}(N+1) \equiv t_N\odot  \left( t_{N-1}\odot
so_{\k_3,\dots,\k_N}(N-1)\right)\equiv iiso(p,q),\quad\! p+q=N-1,
\ee
For example, the Galilean algebra $iiso(N-1)$
appears for $\k_i=(0,0,1,\dots,1)$.
This pattern continues for $\k_1=\k_2=\k_3=0$, etc.

\noindent
{\bf (4)} $\k_N=0$. Here the algebras  have a semidirect sum
structure:
\be
so_{\k_1,\k_2,\dots,\k_{N-1},0}(N+1) \equiv
t'_N\odot  so_{\k_1,\k_2,\dots,\k_{N-1}}(N) \equiv i'so(p,q)  \quad,
\ee
where now $so(p,q)$ acts on $t'_N$  through the contragredient of the
vector representation, hence the notation with a prime.
Of course these algebras are isomorphic to the ones described in (2) as above
(cf. (\ref{cf})).
A pattern similar to that in (3)  occurs for the cases $\k_N=\k_{N-1}=0$, etc.

\noindent
{\bf (5)} $\k_1=\k_N=0$. They have two different semidirect sum
splittings.
The first is:
\be
so_{0,\k_2,\dots,\k_{N-1},0}(N+1)\equiv t_N\odot  ( {t'}_{\! N-1}\odot
so_{\k_2,\dots,\k_{N-1}}(N-1)) \equiv ii'so(p,q) \quad,
\label{circulo}
\ee
where $p+q=N-1$; $so(p,q)$ acts on $t'_{\! N-1}$  through the contragredient
of the vector representation while $i'so(p,q)$ acts on ${t}_{N}$ through
the vector representation. The other is:
\be
so_{0,\k_2,\dots,\k_{N-1},0}(N+1)\equiv t'_N\odot  ( {t}_{\! N-1}\odot
so_{\k_2,\dots,\k_{N-1}}(N-1))  \equiv i'iso(p,q) \quad.
\ee
One example of (\ref{circulo})
is $ii'so(3)$, the Carroll algebra in
$(3+1)$ dimensions \cite{BLL}, which corresponds to $(0,1,1,0)$.

\noindent
{\bf (6)} $\k_a=0$, $a\ne 1,N$.
The structure of these
algebras can be schematically described as $t_r \odot (so(p,q)\oplus
so(p',q'))$  (see \cite{WB}).
In particular, for $\k_2=0$ we have
$t_{2N-2}\odot (so(p,q)\oplus so(p',q'))$ with $p+q=2$ and $p'+q'=N-1$,
which include for $q'=0$ the oscillating  and
expanding  Newton--Hooke algebras \cite{BLL} associated to
$(1,0,1,\dots,1)$ and $(-1,0,1,\dots,1)$, respectively.

\noindent
{\bf (7)} The fully contracted case in the CK family corresponds to
setting all constants  $\k_a=0$. This is the so-called flag algebra
$so_{0,\dots,0}(N+1)\equiv i\dots iso(1)$ \cite{Ros}.

The kinematical algebras associated to different models of space-time
\cite{BLL} belong to the family of CK algebras, and this indeed
provides  one of the strongest physical motivations to study this
family of algebras; in relation with the graded contraction
point of view, see also \cite {MonPatTol} and \cite{HerSanGSGC}.


\sect {Central extensions of the CK algebras}

Our aim in this section is to obtain the general solution to the
problem of finding all central extensions for all the CK algebras and
in arbitrary dimensions.

We write the independent
commutation relations of $so_{\k_1,\dots,\k_N}(N+1)$ (\ref{cb}) as
\be
[\om_{ab},\om_{cd}]
= {\sum_{{\scriptstyle i,j=0\atop \scriptstyle i<j}}^N}
C_{ab,cd}^{ij}\om_{ij}\quad,
\label{da}
\ee
where, as before,
in any $\om_{ef}$ $e<f$ is always assumed
The
four types of structure constants in
(\ref{cb}) are given by
\be
C_{ab,ac}^{ij}=\delta^i_{b}\delta^j_{c}\k_{ab}\quad,\quad
C_{ab,bc}^{ij} = - \delta^i_{a}\delta^j_{c} \quad,\quad
C_{ac,bc}^{ij} =\delta^i_{a}\delta^j_{b}\k_{bc}\quad,\quad a<b<c \quad,
\label{db}
\ee
together with $C_{ab,cd}^{ij}=0$ if all $a,b,c,d$ are different.

Any central extension  $\overline{so}_{\k_1,\dots,\k_N}(N+1)$ of the
algebra $so_{\k_1,\dots,\k_N}(N+1)$ by the one-dimensional algebra of
generator $\Xi$
will have generators $(\om_{ab},\Xi)$ and commutators:
\be
[\om_{ab},\om_{cd}] =
    {\sum_{{\scriptstyle i,j=0\atop \scriptstyle i<j}}^N}
C_{ab,cd}^{ij}\om_{ij}+ \cext_{ab,cd} \Xi\quad,\quad
[\Xi,\om_{ab}]=0\quad,
\label{dcnew}
\ee
where the extension coefficients to be determined
$\cext_{ab,cd}$ (`central charges'), must be antisymmetric in the
interchange of pairs $ab$ and $cd$,
\be
\cext_{cd,ab}=-\cext_{ab,cd} \quad,
\label{ddnew}
\ee
and must fulfil the conditions
\be
{\sum_{{\scriptstyle i,j=0\atop \scriptstyle i<j}}^N}
\left( C_{ab,cd}^{ij}\cext_{ij,ef}+
C_{cd,ef}^{ij}\cext_{ij,ab}+C_{ef,ab}^{ij}\cext_{ij,cd} \right) =0\quad,
\label{de}
\ee
which follow from the Jacobi identity.
The `extension coefficients' are the coordinates
$(\cext(\Omega_{ab},\Omega_{cd})=\cext_{ab,cd})$ of the antisymmetric rank two
tensor $\cext$ which is the two-cocycle of the specific extension
being considered, and  (\ref{de}) is
the two-cocycle condition for the Lie algebra cohomology.
The   classes of
non-trivial two-cocycles associated with the tensors $\alpha$
determine the dimension
of the second cohomology group $H^2(\G,\Re)$.
The rest of this section will be devoted to characterise first the vector space
of all tensors $\cext$ satisfying conditions (\ref{ddnew}) and
(\ref{de}) for the CK algebra $so_{\k_1,\dots,\k_N}(N+1)$.
Secondly, the question of the possible
equivalence of two extensions given by two tensors
$\cext$ will be addressed and solved. The reader who is not
interested in the details of the calculation procedure may skip the
rest of this section; its results are summarised in the theorem at the
beginning of sec. 4.

\subsect{Setting-up the problem: Jacobi identities}

The antisymmetry of $\cext_{ab,cd}$
will be automatically taken into account by considering as independent
coefficients only
those $\cext_{ab,cd}$  with $a \leq c$ and $b<d$ when $a=c$
(the conditions $a<b$ and $c<d$ are always assumed).

The first step consists of solving the ${(N+1)N/2}\choose 3$ (see below)
Jacobi identities (\ref{de}) understood as equations in the coefficients
$\cext_{ab,cd}$.
There is one equation for each possible set of three index
pairs $ab,cd,ef$ in (\ref{de}),
and looking at how many of these indices are {\em
different} we can group all Jacobi identities into four classes:

\noindent
$\bullet$ One equation for the  only possible set $(ab,ac,bc)$ of six indices
made up from  three different indices,
(permutations will lead to the same Jacobi identity).

\noindent
$\bullet$ 16 equations, one for  each possible set of six indices
made up from four different indices
\bea
&&(ab,ac,ad),(ab,ac,bd),(ab,bc,ad),(ab,ac,cd),(ac,ad,bc),(ab,ad,cd),
\nonumber\\
&&(ac,ad,bd),(ab,bc,bd),(ab,bc,cd),(ac,bc,bd),(ab,bd,cd),(ad,bd,bc),
\nonumber\\
&&(ac,bc,cd),(ac,bd,cd),(ad,bc,cd),(ad,bd,cd).
\label{nuevaa}
\eea

\noindent
$\bullet$ 30 equations, one for  each possible set of six indices
made up from five different indices
\bea
&&(ab,ac,de),(ab,ad,ce),(ab,ae,cd),(ac,ae,bd),(ac,ad,be),(ad,ae,bc),\nonumber\\
&&(ab,bc,de),(ab,bd,ce),(ab,be,cd),(ad,bc,be),(ac,bd,be),(ae,bc,bd),\nonumber\\
&&(ab,cd,ce),(ac,bc,de),(ad,bc,ce),(ae,bc,cd),(ac,bd,ce),(ac,be,cd),\nonumber\\
&&(ab,cd,de),(ac,bd,de),(ae,bd,cd),(ad,bc,de),(ad,bd,ce),(ad,be,cd),\nonumber\\
&&(ab,ce,de),(ac,be,de),(ae,bc,de),(ad,be,ce),(ae,be,cd),(ae,ce,bd).
\label{nuevab}
\eea

\noindent
$\bullet$ 15 equations, one for each possible set of six indices
$(ab,cd,e\!f),(ab,ce,d\!f),\dots, $ made up from six different indices.

The following relation, where it is understood that
${m\choose n}=0$ if $m<n$, checks the above splitting:
\be
{{(N+1)N/2}\choose 3} = 1  {{N+1}\choose 3}+16  {{N+1}\choose 4}+
30  {{N+1}\choose 5} +15  {{N+1}\choose 6}.
\label{df}
\ee

Now we write explicitly the above equations.
To begin with, the equation involving only three different indices as
well as the 15 equations with six different indices are easily seen to
be trivially satisfied due to (\ref{ddnew}) or  to the fact that
$C^{ij}_{ab,cd}=0$ whenever $(ab,cd)$ are all different
indices.
The 16 Jacobi
identities for four indices $a<b<c<d$
lead to ten equations, written here in bulk
(later we shall write these equations in a neater way):
\bea
&&\cext_{ab,ad}=-\k_{ab}\cext_{bc,cd} \qquad
\cext_{ac,ad}=\k_{ab}\cext_{bc,bd}  \cr
&&\cext_{ad,cd}=-\k_{cd}\cext_{ab,bc}\qquad
\cext_{ad,bd}=\k_{cd}\cext_{ac,bc} \cr
&&\cext_{ac,cd}=\cext_{ab,bd} \qquad\qquad
\k_{cd}\cext_{ab,ac}=\k_{ab}\cext_{bd,cd} \label{dg}\\
&&\cext_{ac,bd}=-\k_{bc}\cext_{ab,cd}\qquad \cext_{ad,bc}=0\cr
&&\k_{ac}\cext_{ab,cd}=0\qquad\quad\quad
\k_{bd}\cext_{ab,cd}=0\quad.
\nonumber
\eea
(the 16 equations involve several pairs $\k\cext=0$ and $\k\k'\cext=0$; in
these cases the second is clearly a consequence of the first, and may be
therefore
discarded). On their part, the 30 Jacobi identities for five indices
$a<b<c<d<e$, give rise to sixteen equations:
\bea
&&\cext_{ab,ce}=0\cr
&&\cext_{ac,be}=\cext_{ac,de}=0\cr
&&\cext_{ad,be}=\cext_{ad,ce}=0\cr
&&\cext_{ae,bd}=\cext_{ae,bc}=\cext_{ae,cd}=0  \label{dh}\\
&& \k_{bc}\cext_{ab,de}=0\qquad  \k_{cd}\cext_{ab,de}=0\qquad
\k_{de}\cext_{ab,cd}=0 \cr
&&\k_{de}\cext_{ac,bd}=0\qquad \k_{de}\cext_{ad,bc}=0\cr
&&\k_{ab}\cext_{bc,de}=0\qquad \k_{ab}\cext_{bd,ce}=0
\qquad
 \k_{ab}\cext_{be,cd}=0.
\nonumber
\eea

\subsect{Solving strategy}

The structure behind equations (\ref{dg}) and (\ref{dh}) is
not readily apparent. In order to unveil this structure and
solve the central extension problem, we shall

\begin{itemize}

\item[{(1)}] Sort out all extension coefficients (coordinates of $\cext$)
into disjoint classes,

\item[{(2)}] Group all equations as related  to the former
classification and isolate coefficients which can be simply expressed
in terms of the remaining ones,

\item[{(3)}] State the form of the general  solution in terms of
basic extension coefficients from which all others are derived,
but which still may  be subjected to some additional relations, and

\item[{(4)}] Analyse when an extension is trivial (or when two
extensions can be equivalent).

\end{itemize}

It is convenient to put the coordinates of the generic $\cext$ into three
{\it types}:
\smallskip

\noindent
{\it {Type I}}. Coefficients $\cext_{ab,bc}$  with {\it three}
different indices $a<b<c$ with the {\it middle} index common to both
pairs.

\noindent
{\it {Types IIF/IIL}}. Coefficients  $\cext_{ab,ac}$ / $\cext_{ac,bc}$,
with {\it three} different indices $a<b<c$ where
the {\it first / last} index is common to both pairs.

\noindent
{\it {Type III}}. Coefficients  $\cext_{ab,cd}$ with {\it four}
different indices $a<b, \  a<c<d$.

Before studying them below, it may be worthwhile advancing now that,
for all CK algebras, those of type I will always correspond to
two-coboundaries. Those of type IIF/IIL will  determine
two-cocycles (which may be trivial) which result from the
contraction of two-coboundaries through the pseudoextension mechanism
\cite{AldAz,AI}, and those of type  III will determine non-trivial
two-cocycles which are not generated by contraction
from two-coboundaries ({\it i.e.}, different from those in IIF/IIL).

The Jacobi equations (\ref{de}) include a block of
relations like (\ref{dg}) for each set of {\em four} different
indices  $a<b<c<d$, and a group of equations (\ref{dh}) for each set
of {\em five} indices $a<b<c<d<e$.  In order to deal with these
equations simultaneously,  it will be convenient to start with any set
of five indices, and to consider jointly eqs. (\ref{dh}) and five
copies of eqs. (\ref{dg}), one for every choice of four indices out of
the five $abcde$, namely
$abcd$, $abce$, $abde$, $acde$  and $bcde$.  The complete set of
equations thus obtained involves the $45$ central extension
coefficients $\cext_{ab,cd}, \cext_{ab,ce}, \dots$ with indices in the
set $\{a,b,c,d,e\}$ (there are ${5\choose 3}=10$ type I coefficients,
10+10 type IIF+IIL coefficients and 15 (see (\ref{dh})) of type III).
We now write all equations (assuming $a<b<c<d<e$) and group them
in a convenient way.
Out of these 45 coefficients, 30 are either equal
to zero or can be expressed by simple relations in terms of the
remaining extension coefficients.
These 30 coefficients are called $abcde$-{\it derived}, and
are related to the remaining 15 $abcde$-{\it primary} coefficients by the
equations below,
the left/right hand side of which involve only
derived/primary coefficients:
\be
{\rm{I/I}} \qquad
\cext_{ac,cd}= \cext_{ab,bd}  \quad
\cext_{ac,ce}= \cext_{ab,be}  \quad
\cext_{ad,de}= \cext_{ab,be}  \quad
\cext_{bd,de}= \cext_{bc,ce}  \label{dia}
\ee
\bea
{\rm{IIF/IIF}} \qquad &&
\cext_{ac,ad}=\k_{ab}\cext_{bc,bd} \quad
\cext_{ad,ae}=\k_{ac}\cext_{cd,ce} \quad
\cext_{bd,be}=\k_{bc}\cext_{cd,ce} \quad  \label{dib} \\
{\rm{IIL/IIL}} \qquad &&
\cext_{ae,be}=\k_{ce}\cext_{ac,bc} \quad
\cext_{ad,bd}=\k_{cd}\cext_{ac,bc} \quad
\cext_{be,ce}=\k_{de}\cext_{bd,cd} \quad  \label{dic}
\eea
\bea
{\rm{IIF/I}} \qquad &&
\cext_{ab,ad}=-\k_{ab}\cext_{bc,cd} \qquad
\cext_{ab,ae}=-\k_{ab}\cext_{bc,ce} \qquad \cr
\qquad &&
\cext_{ac,ae}=-\k_{ac}\cext_{cd,de} \qquad
\cext_{bc,be}=-\k_{bc}\cext_{cd,de} \qquad  \label{did}
\eea
\bea
{\rm{IIL/I}} \qquad &&
\cext_{ad,cd}=-\k_{cd}\cext_{ab,bc} \qquad
\cext_{ae,ce}=-\k_{ce}\cext_{ab,bc} \qquad \cr
\qquad &&
\cext_{ae,de}=-\k_{de}\cext_{ab,bd} \qquad
\cext_{be,de}=-\k_{de}\cext_{bc,cd} \qquad  \label{dif}
\eea
\be
{\rm{III/III}} \qquad
\cext_{ac,bd}=-\k_{bc}\cext_{ab,cd} \qquad
\cext_{bd,ce}=-\k_{cd}\cext_{bc,de} \qquad  \label{dig}
\ee
\vspace*{0.3cm}
\be
{\rm{III}} \qquad
\begin{array}{ll}
\cext_{ab,ce}=\cext_{ac,be}=\cext_{ae,bc}=0 &
\cext_{ad,bc}=0 \\[0.3cm]
\cext_{ae,bd}=\cext_{ad,be}=0
&
\cext_{be,cd}=0
\\[0.3cm]
\cext_{ac,de}=\cext_{ad,ce}=\cext_{ae,cd}=0 \quad.
\end{array}
\label{dih}
\ee
Sorted out by their type, the 15 $abcde$-primary coefficients are:
\bea
\rm{I}\qquad &&\cext_{ab,bc}
\qquad \cext_{ab,bd} \qquad \cext_{ab,be} \qquad
               \cext_{bc,cd}
\qquad \cext_{bc,ce} \qquad \cext_{cd,de}       \cr
\rm{IIF} \qquad &&\cext_{ab,ac}
 \qquad  \cext_{bc,bd} \qquad \cext_{cd,ce}   \cr
\rm{IIL} \qquad &&\cext_{ac,bc}
 \qquad \cext_{bd,cd} \qquad \cext_{ce,de}    \label{dja} \\
\rm{III}\qquad &&\cext_{ab,cd}
\qquad \cext_{ab,de} \qquad \cext_{bc,de} \nonumber
\eea
These primary coefficients are themselves constrained by the relations:
\be
{\rm{IIF/IIL}} \qquad\qquad
\k_{cd}\cext_{ab,ac}=\k_{ab}\cext_{bd,cd}  \qquad
\k_{de}\cext_{bc,bd}=\k_{bc}\cext_{ce,de} \quad,       \label{djb}
\ee
\bea
&&
\k_{ac}\cext_{ab,cd}=0  \qquad
\k_{bd}\cext_{ab,cd}=0  \qquad
\k_{de}\cext_{ab,cd}=0         \cr
\rm{III} \qquad && \k_{bc}\cext_{ab,de}=0  \qquad
\k_{cd}\cext_{ab,de}=0        \label{djc}\\
&&
\k_{ab}\cext_{bc,de}=0  \qquad
\k_{bd}\cext_{bc,de}=0  \qquad
\k_{ce}\cext_{bc,de}=0 \quad,\nonumber
\eea
which follow from (\ref{dg}), (\ref{dh}). Notice that since the $\k$'s
may be zero, we cannot express any of these primary coefficients in
terms of the others.

We have indicated the type of each group of
coefficients, by the corresponding symbol
in each line, {\it e.g.}  the rows in (\ref{dja})
correspond to the three types  (I, IIF/IIL, III) as given above.
Summing up, the 15 $abcde$-{\it primary} extension coefficients are:

\begin{itemize}

\item[$\bullet$]  Those type I with the first pair $abcde$-contiguous
({\it i.e.}, the indices in the first pair are consecutive in the ordered
sequence $abcde$).

\item[$\bullet$]  Those type IIF/L with the first/last pair of indices
$abcde$-contiguous and three $abcde$-consecutive indices (of course,
the later condition implies the former).

\item[$\bullet$]  Type III with two $abcde$-contiguous pairs.
\end{itemize}

All others are $abcde$-derived.

\subsect{The basic extension coefficients}

The next step in this process is to consider the above
results for all possible numerical values of the five indices $abcde$
since  these numerical values appear as the indices labelling the
coordinates of the tensor $\alpha$. Consider for instance the coordinate
$\alpha_{13,14}$ in a case with, say, $N=8$. This is a $13457$-primary one,
which
means that it cannot be expressed in term of another coefficient with indices
taken from the set $13457$. However, the same coefficient appears as a derived
one for the set of indices $12345$, because the first index pair $13$ are not
contiguous indices in the set $12345$.
It is clear that only those
coefficients with the form of (\ref{dja}) for {\it all} choices of five indices
will be the  really primary, or basic, ones; all other can be ultimately
derived in terms
of these.  We shall call them {\em basic coefficients} of the extension,
and introduce a notation which highlights their
r\^ole in the extensions. By checking the former list of $abcde$-primary
extension coefficients, we readily conclude that

\bigskip
\noindent
{\bf Proposition \propi}.
\\
{\it The {\it basic coefficients} of the extension are:

\noindent
{\it Type I} with the first pair of indices contiguous
\be
\tau_{ac}:= \cext_{a\, a+1,a+1\, c} \qquad a=0,1,\dots,N-2 \quad,
\quad c=a+2,\dots,N \quad, \quad N\ge 2.
\label{maa}
\ee
We remark that the indices in $\tau_{ac}$ cannot be consecutive; there
are
$N(N-1)/2$ basic type I $\tau_{ac}$ extension coefficients.

\noindent
{\it Type  IIF/L}  with three consecutive indices (and
therefore with the first/last pair of contiguous indices):
\bea
&&\cextiif_{a+1\,a+2}:=\cext_{a\, a+1, a\, a+2}\quad,
            \quad a=0,1, \dots, N-2\quad, \quad N\ge 2. \label{maaa}\\
&&\cextiil_{a-2\,a-1}:=\cext_{a-2\, a, a-1\, a}\quad,
            \quad a=2, \dots, N\quad,   \quad N\ge 2. \label{mab}
\eea
There are $(N-1)$ basic type II extension coefficients for each
subtype IIF and IIL.

\noindent
{\it Type III}  with two contiguous pairs of indices:
\be
\cextiii_{b+1\,d+1}:= \cext_{b\, b+1,d\, d+1} \qquad
b=0,1,\dots,N-3\quad,\quad d=b+2,\dots,N-1\quad,
\quad N\ge 3.
\label{mac}
\ee
These $\cextiii$ extension coefficients must have two not consecutive
indices, and the index $0$ cannot appear in any $\cextiii$.
The possible number of these extension coefficients is
$(N-1)(N-2)/2$. }

\bigskip

These basic coefficients are still not independent and must
fulfil the Jacobi relations  (\ref{djb}) and (\ref{djc}).
In particular, for basic type II coefficients, (\ref{djb}) now
reads:
\be
{\rm{IIF/IIL}} \qquad \k_{a+3}\cextiif_{a+1\,a+2} =
\k_{a+1}\cextiil_{a+1\,a+2}\quad, \quad a=0, \dots, N-3\quad,
\label{lfTypeIICondition}
\ee
and (\ref{djc}) for basic Type III coefficients reduces to either
\be
\omega \cextiii_{b+1\,b+3} = 0\quad \hbox{for}\quad
\omega=\k_{b}\ ,\ \k_{b+1}\k_{b+2}\ , \
 \k_{b+2}\k_{b+3}\ ,\ \k_{b+4}\quad,
\label{mbba}
\ee
where for $b=0$/$b$=$N-3$ the first/last condition
$\k_{b}\cextiii=0$/$\k_{b+4}\cextiii=0$ (which would read
$\k_{0}\cextiii=0$/$\k_{N+1}\cextiii=0$) is not present, or to:
\be
\omega \cextiii_{b+1\,d+1}  =0
\quad \hbox{for}\quad
\omega=\k_{b}\ ,\ \k_{b+2}\ ,\ \k_{d}\ ,\ \k_{d+2}
\quad,
\label{mbbb}
\ee
where $b+1$ and $d+1$ are not next neighbours, with similar restrictions as
before for the
extreme equations.

The basic type II coefficients are $\cext_{a\,
a+1, a\, a+2}$ and
$\cext_{a+1\, a+3, a+2\, a+3}$. Both appear in the extended commutators:
\bea
&& [\om_{a\,a+1}, \om_{a\, a+2}] = \k_{a+1} \om_{a+1\, a+2} +
\cext_{a\, a+1, a\, a+2} \Xi \quad,\cr
&& [\om_{a+1\,a+3}, \om_{a+2\, a+3}] = \k_{a+3} \om_{a+1\, a+2} +
\cext_{a+1\, a+3, a+2\, a+3} \Xi \quad,\label{oe}
\eea
so both extension coefficients appear related to the generator
$\om_{a+1\, a+2}$, which explains the notations
$\cextiif_{a+1\,a+2}:=\cext_{a\, a+1, a\, a+2}$ and
$\cextiil_{a+1\,a+2}:=\cext_{a+1\, a+3, a+2\, a+3}$ in (\ref{maaa}),
(\ref{mab}).
Type II basic coefficients are grouped in a  {\it single} coefficient,
$\cextiil_{01}$,  $(N-2)$ {\it pairs},
$\cextiif_{12},\cextiil_{12}$, \dots , $\cextiif_{N-2\,
N-1},\cextiil_{N-2\,N-1}$, and another {\it single} coefficient,
$\cextiif_{N-1\,N}$, the single ones appearing for the cases where the
index pair has not a predecessor or a successor.
Type III basic extension coefficients $\cext_{b\,b+1,
d\,d+1}$ appear in the extended commutators:
\be
[ \om_{b\,b+1}, \om_{d\,d+1} ] =  \cext_{b\,b+1, d\,d+1} \Xi \quad .
\ee

Let us comment on the process of finding the derived coefficients in
terms of the basic extension coefficients. For type I coefficients, the only
constraint is equation (\ref{dia}), which says that all  type I coefficients
$\cext_{ab,bc}$ with the same $ac$ indices are simply equal; this is
$\cext_{ab,bc}=\cext_{a a+1,a+1 c} = \tau_{ac}$. Consider now the
derived type IIF  coefficients.
They must have three non-consecutive indices, and there are three
possibilities, representatives of which
are  the coefficients $\cext_{ab, ad}$ (when
$ab$ are contiguous but $bd$ are not),
$\cext_{ac, ad}$ (when $ac$  are not contiguous but $cd$ are), and
$\cext_{ac, ae}$ (when neither $ac$ nor $ce$ are contiguous).  For the
first one, an equation in (\ref{did}) gives
$\cext_{ab, ad}=-\k_{ab}\cext_{bc, cd}$ and choosing
$c=b+1$ we get $\cext_{ab, ad}=-\k_{ab}\cext_{b\, b+1,b+1\,d}=
-\k_{ab}\tau_{bd}$; as in this case $b=a+1$ we get $\cext_{a\,a+1,
ad}=-\k_{a+1}\tau_{a+1\,d}$ for $d=a+3, \dots, N$.  For the second, one
of equations (\ref{dib}) with $b=c-1$ gives
$\cext_{ac, ad}=\k_{a\,c-1}\cext_{c-1\,c, c-1\,d}$; now as here
$d=c+1$,  and the extension coefficient $\cext_{c-1\,c, c-1\,d}$ is
equal to $\cextiif_{c\,c+1}$, we finally  obtain
$\cext_{ac, a\,c+1}=\k_{a\,c-1}\cextiif_{c\,c+1}$ for $c=a+2, \dots, N-1$.
In the third case, we use one of equations  in (\ref{did}) with
$d=c+1$ to get directly $\cext_{ac, ae} = -\k_{ac}\ \cext_{c\,c+1,
c+1\,e} = -\k_{ac}\tau_{ce}$ for $c=a+2, \dots, N-2$ and $e=c+2, \dots,
N$.
A completely similar process gives the derived type IIL.

Finally, type III has a single class of derived coefficients
which might be different from  zero: $\cext_{ac, bd}$ where
$ab$, $bc$ and $cd$ are contiguous pairs (so all four
$abcd$ indices are consecutive, say $a\,a+1\,a+2\,a+3$). These are
given by equations (\ref{dig}) in term of the basic ones as
$\cext_{a\,a+2, a+1\,a+3}=-\k_{a+1\,a+2}\cext_{a\,a+1, a+2\,a+3} =
-\k_{a+2}\cextiii_{a+1\,a+3}$.  All other non-basic type III
coefficients are necessarily equal to zero.

To sum up, for $N=2$ there are no derived extension coefficients;
for any $N\geq 3$, the complete list of {derived} extension
coefficients is therefore given by

\bigskip

\noindent
{\bf Proposition \propii}.
\\
{\it For $N\geq 3$, the derived extension coefficients are of:

\noindent
{\it Type I}, with the first pair non-contiguous:
\be
\cext_{ac,cd}=\tau_{ad} \quad
a=0, 1, \dots, N-3, \quad
c=a+2, \dots, N-1, \quad  d=c+1, \dots, N.
\label{naa}
\ee

\noindent
{\it Type IIF/L}, with three non-consecutive indices. There are three
possibilities, according to whether the first and second indices or the
second and third are or not contiguous:
\be
\begin{array}{ll}
\cext_{a\,a+1, ad}=-\k_{a+1}\tau_{a+1\,d}&\quad
 a=0,1, \dots, N-3,\quad d=a+3, \dots, N;\cr
\cext_{ac, a\,c+1}=\k_{a\,c-1}\cextiif_{c\,c+1}&\quad
a=0,1, \dots, N-3,\quad c=a+2, \dots, N-1;\cr
\cext_{ac, ae} = -\k_{ac} \tau_{ce}&\quad
a=0,1, \dots, N-4,\quad c=a+2, \dots, N-2,  \cr
&\quad \quad e=c+2, \dots, N \quad  (\hbox{here } N\ge 4).
 \end{array}
\label{nab}
\ee
\be
\begin{array}{ll}
\cext_{a\,c+1, c\,c+1}=-\k_{c+1}\tau_{ac}&\quad
a=0,1, \dots, N-3,\quad c=a+2, \dots, N-1;\cr
\cext_{ac, a+1\,c}=\k_{a+2\,c}\cextiil_{a\,a+1}&\quad
a=0,1, \dots, N-3,\quad c=a+3, \dots, N;\cr
\cext_{ae, ce} = -\k_{ce} \tau_{ac}&\quad
 a=0,1, \dots, N-4,\quad c=a+2, \dots, N-2,  \cr
&\quad \quad e=c+2, \dots, N \quad  (\hbox{here } N\ge 4).
 \end{array}
\label{nac}
\ee

\noindent
{\it Type III}, with at least a non-contiguous pair. Only those of the
form $\cext_{ac,bd}$ with $abcd$ consecutive are possibly different
from zero, and are given by
\be
\cext_{a\, a+2,a+1\, a+3}=-\k_{a+2}\, \cextiii_{a+1\,a+3}
\quad,\quad a=0,1,\dots,N-3 \quad;
\label{nad}
\ee
all other non-basic type III extension coefficients are necessarily
equal to zero.}

It can be checked that for any choice of the extension coefficients
(satisfying the equations (\ref{lfTypeIICondition}), (\ref{mbba}) and
(\ref{mbbb})), the expressions given above for the derived extension
coefficients satisfy all Jacobi equations.
This is cumbersome but straightforward, and will not be done here.

\subsect{Equivalence of extensions: two-coboundaries}

So far we have determined the general form of a two-cocycle on the CK
algebra $so_{\k_1,\dots,\k_N}(N+1)$.  Two two-cocycles differing by a
two-coboundary lead to equivalent extensions,  so the next step is
to find the general form of a coboundary. Let us make the change of
generators  $\om_{ab} \to \om^\prime_{ab} =
\om_{ab}+ \mu_{ab}\Xi$, where $\mu_{ab}$ are arbitrary real numbers.
The commutation relations for the new
generators $\om^\prime_{ab}$, obtained from (\ref{dcnew})  with a given
two-cocycle $\cext_{ab,cd}$ are:
\be
[\om^\prime_{ab},\om^\prime_{cd}]
={\sum_{i,j=0}^N} C_{ab,cd}^{ij}\om^\prime_{ij} + (\cext_{ab,cd} -
{\sum_{i,j=0}^N} C_{ab,cd}^{ij} \mu_{ij}) \Xi\quad.
\label{oa}
\ee
Therefore, the general expression of a two-coboundary $\delta\mu$
generated  by $\mu$ is
\be
(\delta\mu)_{ab,cd} = {\sum_{i,j=0}^N} C_{ab,cd}^{ij} \mu_{ij}   \quad.
\label{ob}
\ee
Using the expressions (\ref{db}) for  the structure constants, we
obtain:
\bea
{\rm I} \qquad  && (\delta\mu)_{ab,bc} = - \mu_{ac} \cr
{\rm IIF/IIL} \qquad && (\delta\mu)_{ab,ac} = \k_{ab}\mu_{bc} \quad/ \quad
(\delta\mu)_{ac,bc} = \k_{bc}\mu_{ab} \cr
 {\rm III} \qquad  && (\delta\mu)_{ab,cd} =
(\delta\mu)_{ac,bd} = (\delta\mu)_{ad,bc} = 0 \quad. \label{oc}
\eea

The question of whether the previously found extension
coefficients (or  two-cocycles) define trivial central extensions
amounts to checking  whether they have the form of a two-coboundary,
(\ref{oc}), which may then be used to eliminate the central $\Xi$ term from
(\ref{oa}).
This depends on the vanishing of the constants
$\k_i$.
In fact, the previous analysis classifies the
extensions into three types, which behave in three
different ways:

\noindent
$\bullet$
Type I extensions can be done for all CK algebras, as
there are no any $\k_i$-dependent restrictions for the basic
type I coefficients $\tau_{ac}$. However these extensions are always
trivial (for all CK algebras simultaneously, as seen in (\ref{oc})),
and will be discarded.
All expressions simplify considerably
if we take this into account, as we shall do from now on.
This `uses up' those coboundaries coming from the values
$\mu_{ac}$ with two non-consecutive
$ac$ indices.
Further equivalences (already for type II) are restricted to
redefinitions of generators with two consecutive indices,
$\om_{a\, a+1} \to
\om^\prime_{a\, a+1} = \om_{a\, a+1} + \mu_{a\, a+1} \Xi$
(see (\ref{og}), (\ref{inI}) below).

\noindent
$\bullet$ {Type II} coefficients can appear  in all CK algebras, as the
$\k_i$-dependent restrictions (\ref{lfTypeIICondition}) are not
strong enough to force all these coefficients to vanish.
However, the triviality of these extensions is also
$\k_i$-dependent, and we will see that the $2(N-1)$
extensions corresponding to the basic extension coefficients
$\cextiil_{01},  \cextiif_{12},\cextiil_{12}, \dots,
\allowbreak
\cextiif_{N-2\,N-1},\cextiil_{N-2\,N-1}, \cextiif_{N-1\,N}$
are all trivial
for  the simple algebras, and all non-trivial for the extreme case of
the flag algebra. It is within this particular type of extensions
that  a {\em pseudoextension} (trivial extension by a
two-coboundary) may become a  non-trivial extension by contraction.

\noindent
$\bullet$ {Type III} coefficients behave in a completely different way.
The $\k_i$-dependent restrictions (\ref{mbba}), (\ref{mbbb})
on type III basic coefficients force many of these coefficients to
vanish (depending on how many constants $\k_i$ are equal to zero).
Those remaining,
once they are present (that is, allowed), are always non-trivial. This
means that there are no type III non-trivial central extensions coming
by  contraction from pseudoextensions.

All that remains is to discuss the possible equivalence among type II
extensions. The basic type II values of the  coboundary associated to
the change of generators
$\om_{a+1\, a+2} \to
\om^\prime_{a+1\, a+2} = \om_{a+1\,a+2} + \mu_{a+1\,
a+2} \Xi$ are:
\be
{\rm IIF/IIL} \qquad
   (\delta\mu)^{F}_{a+1\,a+2} = \k_{a+1}\mu_{a+1\,a+2} \qquad
   (\delta\mu)^{L}_{a+1\,a+2} = \k_{a+3}\mu_{a+1\,a+2} \quad.
\label{og}
\ee
We remark that (as it should), these coboundaries automatically
satisfy  equation (\ref{lfTypeIICondition}).  We must study now how the
freedom afforded by these changes can be used to reduce to zero some of
the extension coefficients.

Consider first the single $\cextiil_{01}$, the value of which can
be arbitrary. We see in (\ref{og}) that as long as
$\k_{2} \neq 0$, we can reduce it  to zero by using the coboundary
$\mu_{01}$. Then the extension corresponding to the basic coefficient
$\cextiil_{01}$ is non-trivial when $\k_{2} = 0$ and trivial otherwise.
Likewise,  the extension corresponding to the basic coefficient
$\cextiif_{N-1\,N}$ is non-trivial when $\k_{N-1} = 0$ and trivial
otherwise.  The possible triviality of these extensions is thus
completely governed by two constants $\k_{i}$ which play a special
r\^ole:  the second $\k_2$ and the last but one $\k_{N-1}$.

Let us now look at the case of pairs
$\cextiif_{a+1\,a+2}, \cextiil_{a+1\,a+2}$. Here the situation is
controlled by two constants
$\k_{i}$, namely $\k_{a+1}$ and $\k_{a+3}$.  When they are both equal
to zero, equation (\ref{lfTypeIICondition}) is automatically satisfied,
irrespective of the values of the pair of coefficients
$\cextiif_{a+1\,a+2}$ and $\cextiil_{a+1\,a+2}$  which cannot be
modified by adding a coboundary; in this case the cocycles
associated to these coefficients are simultaneously non-trivial. When
only one of the constants $\k_{a+1}$ and $\k_{a+3}$ is equal to zero,
equation (\ref{lfTypeIICondition}) forces the vanishing of one of the
coefficients in the pair, while the other can be reduced to zero by the
appropriate coboundary coming from $\mu_{a+1\,a+2}$. And finally, when
both constants $\k_{a+1}$ and $\k_{a+3}$ are different  from zero, then
(\ref{lfTypeIICondition}) enforces the possibility  of simultaneously
reducing $\cextiif_{a+1\,a+2}$ and $\cextiil_{a+1\,a+2}$ to zero by
using the coboundary coming from $\mu_{a+1\,a+2}$. Therefore  the two
two-cocycles extensions corresponding to the two basic coefficients
$\cextiif_{a+1\,a+2}$ and $\cextiil_{a+1\,a+2}$ are non-trivial when
both $\k_{a+1} = 0$ and  $\k_{a+3} = 0$; the two extensions are
simultaneously trivial otherwise.

Once the coboundary Type I coefficients are removed, the contents
of Prop. \propi\ and \propii\ may be summarised by the following
\goodbreak

\noindent
\begin{tabular}{|l|l|}
\hline
\multicolumn{2}{|c|}{\it \hfill Basic coefficients and relations\hfill
Number of them}\\
\hline
{Type IIF/L} & $\cextiif_{a+1\,a+2} :=
\cext_{a\,a+1,a\,a+2}\quad\cextiil_{a\,a+1}
:=
\cext_{a\,a+2,a+1\,a+2} \hfill 2(N-1)$\\
& \quad $a=0,1,\dots,N-2\quad N\ge 2$\\
&$\k_{a+3}\cextiif_{a+1\,a+2} = \k_{a+1}\cextiil_{a+1\,a+2}  \quad
a=0, \dots N-3 $\\
\hline
{Type III} &$\cextiii_{b+1\,d+1} := \cext_{b\,b+1, d\,d+1}
\hfill (N-1)(N-2)/2$\\
& $  \quad b=0,1,\dots,N-3 \quad \quad  d= b+2,\dots, N-1   \quad N\ge
3$\\
& $ \k_{b}\,\cextiii =
\k_{b+1}\k_{b+2}\,\cextiii =
\k_{b+2}\k_{b+3}\,\cextiii =
\k_{b+4}\,\cextiii=0$\ \, for \ \, $\cextiii\equiv\cextiii_{b+1\,b+3}$\\
& $\k_{b}\,\cextiii=\k_{b+2}\,\cextiii =
 \k_{d}\,\cextiii=\k_{d+2}\,\cextiii=0$ \  for \
 $\cextiii\equiv\cextiii_{b+1\,d+1}$\\
 &\hskip 5cm with\ \, $d=b+3,\dots,N-1$\\
\hline
\hline
\multicolumn{2}{|c|}{\it Derived coefficients}\\
\hline
{Type IIF/L}&
$\cext_{ac, a\,c+1}=\k_{a\,c-1}\cextiif_{c\,c+1}$\\
&  \quad $ a=0,1, \dots, N-3 \quad \quad   c=a+2, \dots, N-1   \quad \quad N\ge
3$\\ &$\cext_{ac, a+1\,c}=\k_{a+2\,c}\cextiil_{a\,a+1}$\\
& \quad $ a=0,1, \dots, N-3 \quad \quad c=a+3, \dots, N  \quad \quad N\ge 3$\\
\hline
{Type III} &
$\cext_{a\, a+2,a+1\, a+3}=-\k_{a+2}\, \cextiii_{a+1\,a+3}
\quad a=0,1,\dots,N-3\quad N\ge 3$\\
\hline
\end{tabular}

\smallskip
\noindent
\centerline{
{\bf{Table I}}. Basic and derived Types II and III extension
coefficients for CK algebras}

\bigskip


\sect {The second cohomology groups of the CK algebras}

\subsect{The structure of the central extensions of a CK algebra}

If we completely disregard the type I extensions, which
are trivial for all the CK algebras, all the results obtained
in sec.~3 can be summed up in the following theorem, which
contains the complete solution to the problem of finding the central
extensions of CK algebras:

\bigskip

\noindent
{\bf Theorem \teor}.
\\
{\it The independent non-zero commutators
of any central extension $\overline{so}_{\k_1,\dots,\k_N}(N+1)$ of the
CK Lie algebra
${so}_{\k_1,\dots,\k_N}(N+1)$ can be written as:
\be
\begin{array}{ll}
[\om_{ab}, \om_{bc}]   = - \om_{ac}&\quad \cr
[\om_{ab}, \om_{a\,b+1}]
=\k_{ab} \om_{b\,b+1}+\k_{a\,b-1} \cextiif_{b\,b+1}\Xi &\quad
[\om_{ab}, \om_{ac}]  = \k_{ab} \om_{bc}   \quad \mbox{\rm for}\
c > b+1 \cr
[\om_{ac}, \om_{a+1\,c}]   =\k_{a+1\,c} \om_{a\,a+1}+\k_{a+2\,c}
\cextiil_{a\,a+1}\Xi &\quad
[\om_{ac}, \om_{bc}]  = \k_{bc}\om_{ab} \quad \mbox{\rm for}\ b > a+1  \cr
[\om_{a\,a+1}, \om_{c\,c+1}]  =\cextiii_{a+1\,c+1} \Xi
&\quad[\om_{a\,a+2}, \om_{a+1\,a+3}]  =-\k_{a+2} \cextiii_{a+1\,a+3} \Xi
\quad
\end{array}
\label{pf}
\ee
where $\k_{aa}:=1$. The extension is completely described
by a number of extension coefficients:

\noindent
 $\bullet$   A {\it single} type II coefficient, $\cextiil_{01}$, which
produces an extension which is non-trivial if $\k_2=0$ and trivial
otherwise.

\noindent
 $\bullet$   $(N-2)$ type II {\it pairs},
$\cextiif_{12},\cextiil_{12}$;
\dots ;
$\cextiif_{N-2\, N-1},\cextiil_{N-2\,N-1}$. Each pair of coefficients
must satisfy
$\k_{a+3}\cextiif_{a+1\,a+2} =
\k_{a+1}\cextiil_{a+1\,a+2}$. The two extensions corresponding to the
pair
$\cextiif_{a+1\,a+2}$ and $\cextiil_{a+1\,a+2}$ are both non-trivial
when $\k_{a+1} = 0$ and  $\k_{a+3} = 0$. The two two-cocycles are
simultaneously trivial otherwise.

\noindent
 $\bullet$ A {\it single} type II coefficient, $\cextiif_{N-1\,N}$, which
produces an extension which is non-trivial if $\k_{N-1}=0$ and trivial
otherwise.

\noindent
$\bullet$  $(N-2)$ type III extension coefficients $\cextiii_{13},
\cextiii_{24},
\dots, \cextiii_{N-2\,N}$, satisfying:
\be
\omega \cextiii_{b+1\,b+3} = 0
\quad\hbox{for}\quad
\omega=\k_{b},\; \k_{b+1}\k_{b+2},\; \k_{b+2}\k_{b+3}, \; \k_{b+4},
\label{aI}
\ee
where when either $b=0$ or $b=N-3$ the first or last
condition, which would read
$\k_{0}\cextiii=0$ or $\k_{N+1}\cextiii=0$
is not present. The extension
corresponding to any of these non-zero coefficients is always non-trivial.

\noindent
$\bullet$  $(N-2)(N-3)/2$ type III extension coefficients
$\cextiii_{14}, \cextiii_{15}, \dots, \cextiii_{1N}$;
$\cextiii_{25}, \dots, \cextiii_{2N}$;
$\dots,$ $\dots; \cextiii_{N-3\,N}$
whose indices differ by more than two. The coefficient
$\cextiii_{b+1\,d+1}$ satisfies
\be
\omega\cextiii_{b+1\,d+1}  =0\quad \hbox{for}\quad
\omega=\k_{b},\; \k_{b+2},\; \k_{d},\; \k_{d+2}
\label{aII}
\ee
with similar restrictions as to the actual presence of the
equations involving the non-existent values
$\k_{0}$ or $\k_{N+1}$.
The extension corresponding to any of these non-zero coefficients is
always non-trivial.}

\medskip
All type II extensions come from the
pseudocohomology mechanism \cite{AldAz,AI};
if $\omega_{b+1}\ne 0\ne \omega_{b+3}$ we obtain from
(\ref{pf})
\bea
&&[\Omega_{a\,b+1},\Omega_{a\,b+2}]=\omega_{a\,b+1}
(\Omega_{b+1\,b+2}+{1\over\omega_{b+1}}
\cextiif_{b+1\,b+2}\Xi)
\quad,\cr
&&[\Omega_{b+1\,c},\Omega_{b+2\,c}]=
\omega_{b+2\,c}(\Omega_{b+1\,b+2}+{1\over \omega_{b+3}}
\cextiil_{b+1\,b+2}\Xi)
\label{inI}
\eea
so we may remove  the extension coefficients by means of a redefinition
of the generator  $\Omega_{b+1b+2}$ given by the one-cochain of non-zero
coordinates   ${1\over \omega_{b+1}} \cextiif_{b+1 b+2}$
$= {1\over\omega_{b+3}} \cextiil_{b+1 b+2}$.  When $\omega_{b+1}\
(\omega_{b+3})$  goes to zero the extension given by
$\cextiif_{b+1 b+2}$, $(\cextiil_{b+1 b+2})$  might be non-trivial
(because the one-cochain from it comes diverges but
$\omega_{a\,b+1}/\omega_{b+1}$ does not).  Due to the Jacobi identity
(represented here by the equation  $\omega_{b+3}\cextiif_{b+1\,
b+2}=\omega_{b+1}\cextiil_{b+1\, b+2}$),  the non-trivial extension exists
if both $\omega_{b+1}=\omega_{b+3}=0$.  In terms of the standard
triangular arrangement of generators,  it is worth remarking that each
Type II coefficient $\cextiif_{b+1\,b+2}$  appears only in the extended
commutators of two horizontal neighbours in the  columns of
$\Omega_{0 b+1},\Omega_{0\,b+2}$, while each  Type II
coefficient $\cextiil_{b+1\,b+2}$ appears only in the extended commutators
of two vertical neighbours in the rows of
$\Omega_{b+1 N},\Omega_{b+2 N}$.  The corresponding extension is
non-trivial only when both
$\omega_{b+1}=\omega_{b+3}=0$; this means that the algebra has two
different semidirect structures (cf. sec. 2). This is exhibited by the
two rectangle boxes in the following diagram, where we have shortened
$\cextiif \equiv \cextiif_{b+1\,b+2}$ and $\cextiil \equiv \cextiil_{b+1
b+2}$ and we have also indicated the
$\k_i$ factors which appear in these extended commutators for the
generators which are inside one of the boxes but outside the other:

 \bigskip
\noindent
\begin{tabular}{ccc|ccc|ccccl}
$\om_{01} $& $\ldots$&
$\om_{0\, b} $&
\multicolumn{3}{|c|}{
 $\om_{0\,b+1}$
{${\scriptstyle{\k_{0b}\cextiif}\atop {\displaystyle{\longrightarrow}}}$}
 $\om_{0\, b+2}$}&$\om_{0\, b+3}$&$\om_{0\, b+4}$&
$\ldots$&$\om_{0N}$&\\
&$\ddots $&$\vdots$&
\multicolumn{3}{|c|}{$\vdots\qquad\qquad\qquad\vdots$}&$\vdots$&
$\vdots $&&$\vdots$\\
   &$ $&$\om_{b-1\,b}$&
\multicolumn{3}{|c|}{ $\om_{b-1\,b+1}
 {\scriptstyle{\k_{b-1\,b}\cextiif}\atop {\displaystyle{\longrightarrow} } }
  \om_{b-1\,b+2} $}&$\om_{b-1\,b+3}$&
$\om_{b-1\,b+4}$&$\ldots$&$\om_{b-1\,N}$\\
  &$ $&$ $&\multicolumn{3}{|c|}{
  $\om_{b\,b+1}$
{${\scriptstyle{ \cextiif}\atop {\displaystyle{\longrightarrow}}}$}
 $\om_{b\,b+2}$}&$\om_{b\,b+3}$&$\om_{b\,b+4}$&
$\ldots$&$\om_{b\,N}$&\\
\cline{4-10}
 &$ $&\multicolumn{4}{r|}{$\om_{b+1\,b+2}$}&$\om_{b+1\,b+3}$&
$\om_{b+1\,b+4}$&$\ldots$&
 $\om_{b+1\,N}$  &\\
  &$ $& \multicolumn{3}{c}{\,}   & & $ \downarrow
{ \scriptstyle{ \cextiil } }$&
\multicolumn{2}{c}{$ \downarrow
{ \scriptstyle{\k_{b+3\,b+4}\cextiil } }$}    &
   \multicolumn{2}{r}{ $\quad \downarrow
{ \scriptstyle{\k_{b+3\,N}\cextiil } }$}  \\
  &$ $& \multicolumn{3}{c}{\,} &   &$\om_{b+2\,b+3}$&$\om_{b+2\,b+4}$&
$\ldots$&$\om_{b+2\,N}$&\\
\cline{7-10}
   & \multicolumn{6}{c}{\,}     &$\om_{b+3\,b+4}$ &
$\ldots$&$\om_{b+3\,N}$&\\
   &$ $& \multicolumn{5}{c}{\,}    &  &$\ddots$&$\vdots$\\
   &$ $& \multicolumn{5}{c}{\,}     &  &
$ $&$\om_{N-1\,N}$&\\
\end{tabular}

\bigskip

\subsect{The dimension of the second cohomology groups of the CK contracted
algebras}

Theorem \teor\ contains all
the necessary information to determine easily the
dimension of the second cohomology group
$H^2(so_{\k_1,\dots,\k_N}(N+1), \Re)$ of the CK algebras.
This dimension is obtained
as the sum of a number of completely independent contributions, each
one related to the vanishing of specific groups of constants $\k_i$ as follows:

\noindent
$\bullet$ 1 when $\k_{2} = 0$, with two-cocycle determined by the basic
coefficient $\cextiil_{01}$.

\noindent
$\bullet$  2 for each pair of next neighbour zero constants,
$\k_{1} = \k_{3} = 0$, $\k_{2} = \k_{4} = 0$,  \dots, $\k_{N-2} =
\k_{N}=0$. The two-cocycles appearing with the vanishing pair
$\k_{a+1} = \k_{a+3} = 0$ are determined by basic extension coefficients
$\cextiif_{a+1\,a+2}$ and $\cextiil_{a+1\,a+2}$. This might amount to
a subtotal of $2(N-2)$ when all pairs of second neighbours are zero,
{\it i.e.}, when all constants
$\k_{i}$ are zero.

\noindent
$\bullet$  1 when $\k_{N-1} = 0$, with extension coefficient
$\cextiif_{N-1\,N}$.

\noindent
$\bullet$    1 for each group of zero constants
$\{\k_{b}, \k_{b+2}, \k_{b+4}\}$ or $\{\k_{b}, \k_{b+1}, \k_{b+3},
\k_{b+4}\}$ with extension coefficient $\cextiii_{b+1\,b+3}$.

\noindent
$\bullet$    1 for each group of zero constants
$\{\k_{b}, \k_{b+2}, \k_{d},  \k_{d+2}\}$  with extension coefficient
$\cextiii_{b+1\,d+1}$ for
$d=b+3,\dots ,N-1$.

As mentioned after (\ref{aI}), (\ref{aII}),
the literal application of the two last
rules may apparently involve the constants
$\k_0$, $\k_N$. In these cases the corresponding conditions involving these
inexistent values should be disregarded.

We can translate the previous rules into a closed formula for the dimension of
the second cohomology group
$H^2(so_{\omega_1\dots\omega_N}(N+1),\Re)$. Let $\delta_i$ be defined by
\begin{equation}
\delta_i=\left\{
\begin{array}{l@{\quad}l}
1 & \omega_i=0 \\
0 & \omega_i\ne 0
\end{array}
\right. \quad (i=1, \dots, N),
\label{insert1}
\end{equation}
then dim$H^2(so_{\omega_1\dots\omega_N}(N+1),\Re)$ is given in terms
of the sequence $\delta_1, \delta_2, \dots, \delta_N$ by
\begin{equation}
\begin{array}{l}
\displaystyle
\mbox{dim}(H^2(so_{\omega_1\dots\omega_N}(N+1),\Re))
= \delta_{2} + \delta_{_{N-1}} + 2\sum_{i=1}^{N-2}
\delta_{{i}}\delta_{{i+2}}
\\ \displaystyle \hfill
+\sum_{i=1}^{N-2} \delta_{{i}} \delta_{{i+4}}
[ \delta_{{i+2}} + \delta_{{i+1}} \delta_{{i+3}} -
\delta_{{i+2}} \delta_{{i+1}} \delta_{{i+3}} ]
+ \sum_{i=1}^{N-3} \sum_{j=i+3}^{N}
\delta_{{i}} \delta_{{i+2}}
\delta_{{j}}\delta_{{j+2}}
\end{array}
\label{insert2}
\end{equation}
where $\omega_i\equiv 0$ ($\delta_i = 1$) for $i>N$.
For instance, if $\omega_i=0\ \forall\ i=1,\dots,N$ (flag algebra)
then all $\delta_i=1$, hence all terms in
(\ref{insert2}) contribute and we obtain
\begin{equation}
\begin{array}{l}
\displaystyle
\mbox{dim}(H^2(so_{0\dots 0}(N+1),\Re)) =
2 + 2(N-2) + (N-2) + \sum_{i=1}^{N-3} (N-i-2) =
\\ \displaystyle \hfill
= 2(N-1) + {(N-2)(N-1)\over 2}
={N(N+1)\over 2} -1 \quad.
\end{array}
\label{insert3}
\end{equation}
Each term in formula (\ref{insert2}) is related to a given extension
coefficient as stated in Theorem \teor\ and the preceding rules.

To effectively apply the above rules, it is
convenient to browse
through the list of possible extension coefficients and to see whether
each of them is  allowed/trivial for the algebra we are dealing with or not.
As a first example, consider the algebra  $so_{0,\k_2,0,0}(5)$
with $\k_2\ne 0$. For any CK algebra
$so_{\k_1,\k_2,\k_3,\k_4}(5)$, the possible extension
coefficients are: $\cextiil_{01}$, $\cextiif_{12}$, $\cextiil_{12}$,
$\cextiif_{23}$, $\cextiil_{23}$, $\cextiif_{34}$; $\cextiii_{13}$,
$\cextiii_{14}$, $\cextiii_{24}$. In this case,  it is clear that the
type II non-trivial ones are only $\cextiif_{12}$, $\cextiil_{12}$ (as
$\k_{1} = \k_{3} = 0$) and $\cextiif_{34}$ (as here
$\k_{N-1} \equiv \k_{3}= 0$). Type III extension coefficient
$\cextiii_{13}$ is allowed, and therefore gives a non-trivial cocycle,
as $\k_{1}=0, \k_{3}=0$ and $\k_4=0$. Type III extension
$\cextiii_{14}$ is not allowed, since $\k_{2}$ and $\k_{3}$  are not
simultaneously equal to zero. Type III coefficient $\cextiii_{24}$ is
allowed (and therefore  non-trivial), as  $\k_1=0,\k_3=0$.  So the
dimension of the second cohomology group  is equal to 5 in this case.

The dimension of  $H^2$ for many other algebras can be
derived from these rules.  Although in the next section we shall
give explicitly all the extended CK algebras up to $N=4$, we
mention here the result for some interesting algebras (see sec.
2). In some cases these cohomology groups have been known since
long.

\noindent
{(1)} When all $\k_i$ are different from
zero,  all type III coefficients are equal to zero, and
all type II (which can be different from zero) are coboundaries.
Therefore, the second cohomology group of the simple
pseudo-orthogonal algebra $so_{\k_1,\dots,\k_N}(N+1) \ \forall\,  \k_i \neq 0
$, is
trivial
in accordance with the Whitehead lemma, $H^2(\G,\Re)=0$ if $\G$ is semisimple.

\noindent {(2)} If only the first constant is equal to zero,
$\k_1=0$, we see that the inhomogeneous algebras $iso(p,q)$ where
$p+q=N$ (e.g., Euclidean and Poincar\'e) have no  non-trivial
extensions except in the case
$N=2$, where the first constant
also plays the r\^ole of the last but
one, and there is a single extension coefficient $\cextiif_{12}$. This
result is just a rephrasing of the statement that
in this case every local exponent
is equivalent to zero for $N>2$, as found by Bargmann \cite{Bar} in his
classical study.

\noindent {(3)} When $\k_1=\k_2=0$ (all others being non-zero)
the twice inhomogeneous algebras
$iiso(p,q)$  have a non-trivial extension coefficient: $\cextiil_{01}$
(this is just the mass for the Galilei algebra, which parametrises its
second  cohomology group).  Generically, this is the only non-trivial
extension in this family of algebras, though in the lower
dimensional cases
$N=2, 3$ additional non-trivial extensions appear, as seen in the
examples below.

\noindent {(4)} The flag algebra $ii\ldots iso(1)$ is the most contracted
algebra  in the
CK family, and corresponds to all $\k_i=0$. In this case, basic type II
or III extension coefficients, whenever present, lead to non-trivial
extensions. Furthermore, all the conditions that these coefficients must
satisfy are automatically fulfilled, as a consequence of the vanishing of all
$\k_i$.
There are $2(N-1)$ Type II and $[(N-1)(N-2)/2]$ Type III coefficients
in this case, so that:
\be
H^2(so_{0,\dots,0}(N+1), \Re)= \Re^{2(N-1) + [(N-1)(N-2)/2]}=
\Re^{[{\left(N+1\atop 2\right)}-1]}\quad.
\label{odb}
\ee
The dimension of $H^2$ for the flag algebra is just equal to
$\dim(so_{\k_1,\dots,\k_N}(N+1))-1$ (see (\ref{insert3})).

To conclude this section we mention that, had we considered graded
contractions from $so(N+1)$ beyond the CK family,  we would have found a
larger set of algebras with the [$(N+1)N/2$]-dimensional Abelian
algebra as the most contracted one. Since for it all equations (\ref{de}) are
trivially
satisfied and only the antisymmetry conditions (\ref{ddnew}) remain, the
cohomology group
of this Abelian algebra has dimension
$ \left( N+1\atop 2 \right)  \left( \left( N+1\atop 2 \right) -1
\right) /2$.


\sect {Examples: all central extensions for $N=2,3,4$}

We extract from the general solution in Theorem \teor\
the central extensions for all the
CK algebras $so_{\k_1,\dots,\k_N}(N+1)$ for $N=2,3,4$ \cite{Tesis}.
We remark that our results cover in a single stroke a large family of Lie
algebras; in particular, the family ${so}_{\k_1,\k_2,\k_3,\k_4}(4+1)$ contains
all relativistic and non-relativistic $3+1$ kinematical algebras,
the cohomology of which can be then read off directly.

\subsect{$\overline{so}_{\k_1,\k_2}(3)$}

There are two central extension coefficients of type II:
\be
\cextiil_{01} \equiv
\cext_{02,12}\quad, \quad
\cextiif_{12} \equiv \cext_{01,02}\quad.
\label{fb}
\ee
which are not constrained by any additional condition. The
Lie brackets of  $\overline{{so}}_{\k_1,\k_2}(3)$ are
\be
[\om_{01},\om_{02}] =\k_1 \om_{12}+\cextiif_{12}\Xi \quad, \quad
[\om_{01},\om_{12}] =-\om_{02} \quad,\quad
[\om_{02},\om_{12}] =\k_2 \om_{01}+\cextiil_{01}\Xi \quad.
\label{fc}
\ee
The triviality of any such
extension  is governed by the second and last but one constants in the
list $\k_i$. In this case, these are the second and the first. Thus
$\cextiil_{01}$ is trivial if $\k_2\ne 0$, and
$\cextiif_{12}$ is trivial if $\k_1\ne 0$.  This is exhibited by the
redefinitions:
\be
\om_{01}\to \om_{01}+ \frac{\cextiil_{01} }{\k_2}\Xi\quad,\quad \k_2\ne 0
\quad;\quad
\om_{12} \to \om_{12}+\frac{\cextiif_{12}}{\k_1}\Xi\quad,\quad \k_1\ne 0\quad.
\label{fd}
\ee
Thus, dim[$H^2({so}_{\k_1,\k_2}(3),\Re)$] is equal to:

\noindent
$\bullet$ 0 for the simple Lie algebras $so(3)$ and $so(2,1)$ (both
$\k_1$ and $\k_2\ne 0$).

\noindent
$\bullet$ 1 for the 2-dimensional Euclidean algebra, which appears for
the constants $(0, 1)$ (extension coefficient $\cextiif_{12}$) and
$(1,0)$ (extension coefficient $\cextiil_{01}$).

\noindent
$\bullet$ 1 for the (1+1)-Poincar\'e algebra, which also appears
twice, as  $(0, -1)$ and $(-1,0)$, respectively with extension
coefficients $\cextiif_{12}$ and $\cextiil_{01}$.

\noindent
$\bullet$ 2 for the (1+1)-Galilei algebra, which appear for constants
$(0, 0)$, with both extension coefficients $\cextiif_{12}$ and
$\cextiil_{01}$. Physically, these extensions are
parametrised by the mass and a constant force (see \cite{LL} and
references therein and \cite{AldAzcTuc}).

\subsect{$\overline{so}_{\k_1,\k_2,\k_3}(4)$}

The full set of extension possibilities appears first in this
case.
However, there are some non-generic coincidences. There are
four basic extension coefficients of type II, and one of type III:
\be
 \cextiil_{01}\equiv \cext_{02,12}\quad,\quad
\cextiif_{12}\equiv \cext_{01,02} \quad,\quad
 \cextiil_{12}\equiv \cext_{13,23}\quad,\quad
\cextiif_{23}\equiv \cext_{12,13}\quad,\quad
\cextiii_{13}\equiv  \cext_{01,23} \ ;
\label{fe}
\ee
which must satisfy
\be
\k_3\cextiif_{12}=\k_1\cextiil_{12}\quad,\quad
\k_1\k_2\cextiii_{13}=0\quad, \quad
\k_2\k_3\cextiii_{13}=0\quad.
\label{ff}
\ee
Then, eqs. (\ref{pf}) give the commutation  rules of
$\overline{so}_{\k_1,\k_2,\k_3}(4)$:
\be
\begin{array}{lll}
[\om_{01}, \om_{02}]  =\k_1 \om_{12}+\cextiif_{12}\Xi     &
 \ [\om_{01}, \om_{12}]   = - \om_{02}    &\
 [\om_{02}, \om_{12}]   =\k_2 \om_{01}+\cextiil_{01}\Xi   \cr
 [\om_{01}, \om_{03}]  = \k_{1} \om_{13}   &
 \ [\om_{01}, \om_{13}]  = - \om_{03}   &\
 [\om_{03}, \om_{13}]  = \k_3(\k_2 \om_{01}+\cextiil_{01}\Xi )   \cr
  [\om_{02}, \om_{03}]  =\k_1(\k_2 \om_{23}+\cextiif_{23}\Xi )   &
 \ [\om_{02}, \om_{23}] = -  \om_{03}   &\
 [\om_{03}, \om_{23}]  = \k_{3} \om_{02}  \cr
 [\om_{12}, \om_{13}]  = \k_2 \om_{23}+\cextiif_{23}\Xi    &
 \ [\om_{12}, \om_{23}]  = -  \om_{13}   & \
 [\om_{13}, \om_{23}]  = \k_3 \om_{12}+\cextiil_{12}\Xi   \cr
 [\om_{01},\om_{23}] =\cextiii_{13}\Xi &
 \ [\om_{02},\om_{13}]=-\k_2\cextiii_{13}\Xi  &\
[\om_{03},\om_{12}]=0\quad.\cr
\end{array}
\label{fhh}
\ee

The extension coefficient
$\cextiil_{01}$ produces a  non-trivial cocycle when the second
constant $\k_2 = 0$, the extension $\cextiif_{23}$ is non-trivial when
the last but one constant ($\k_2$ again in this case) is
zero, and the extensions given by $\cextiif_{12}$ and $\cextiil_{12}$ are
non-trivial when $\k_1 =\k_3 = 0$. The extension  determined by
$\cextiii_{13}$ is
only present (see (\ref{ff}))
when $\k_2 = 0$ or when $\k_1 = \k_3 = 0$, and whenever it
appears, it is non-trivial.
The redefinition of  generators displaying the triviality of type II
extensions is:
\bea
&& \om_{01}\to \om_{01} +\frac{\cextiil_{01}}{\k_2}\Xi\quad \mbox{if}\quad
\k_2\ne 0
\cr
 && \om_{12}\to \om_{12} +\frac{\cextiil_{12}}{\k_3}\Xi\quad
\mbox{if}\quad \k_3\ne 0  \qquad
\om_{12}\to \om_{12} +\frac{\cextiif_{12}}{\k_1}\Xi\quad \mbox{if}\quad
\k_1\ne 0 \cr
 && \om_{23}\to \om_{23}+\frac{\cextiif_{23}}{\k_2}\Xi\quad
\mbox{if}\quad \k_2\ne 0 \qquad.
\label{fi}
\eea
Notice that when $\k_1$ and $\k_3$ are  both different from zero,
equation (\ref{ff}) guarantees that both expressions for the
redefinition of $\om_{12}$ indeed coincide.

To conclude the analysis we now give
dim[$H^2(so_{\k_1,\k_2,\k_3}(4),\Re)$]:

\noindent
$\bullet$   3 for $\k_2=0$ with either $\k_1$ or $\k_3$ non-zero:
non-trivial extension coefficients $\cextiil_{01}$,
$\cextiif_{23}$ and
$\cextiii_{13}$. Examples here are both
(2+1) Newton--Hooke algebras $(\pm 1,0,1)$,  $(1,0,\pm 1)$ and the (2+1)
Galilean one
$iiso(2)$   $(0,0,1)$, $(1,0,0)$. Also
$iiso(1,1)$ $(-1,0,0)$, $(0,0,-1)$  and $t_4 \odot (so(1,1)\oplus so(1,1))$
$(-1,0,-1)$.

\noindent
$\bullet$  3 for $\k_1=\k_3=0$ and  $\k_2\ne 0$; non-trivial extensions are
$\cextiif_{12}$,
$\cextiil_{12}$ and $\cextiii_{13}$. Here we find the (2+1) Carroll
algebra $(0,1,0)$ and ${ii'so}(1,1)$ $(0,-1,0)$.

\noindent
$\bullet$ 5 for the most contracted algebra in the CK family with
$\k_1=\k_2=\k_3=0$; it corresponds to the flag space algebra
${iiiso}(1)$.

\noindent
$\bullet$ 0 for all the remaining algebras. Up to isomorphisms these
include the semisimple ones $so(4)$ (once), $so(3,1)$ (four times),
$so(2,2)$ (three times), the three-dimensional Euclidean $iso(3)$ (two times)
and
Poincar\'e  algebras $iso(2,1)$ (six times).

For convenience we include below  the standard triangular diagram with all the
extended
commutators for $\overline{so}_{\k_1,\k_2,\k_3}(4)$.
An arrow  between generators
$A$ and $B$ means that a central $\Xi$-term, with coefficient indicated
near the arrow, is added to the non-extended commutator $[A,B]$. In the
usual kinematical interpretation, the generators may be translated as
$\om_{01} \to H$ (Hamiltonian),
$\om_{02} \to P_1, \
\om_{03} \to P_2$ (momenta),
$\om_{12} \to K_1, \
\om_{13} \to K_2$ (boosts),
$\om_{23} \to J$  (rotation).

\hskip 1.5truecm 
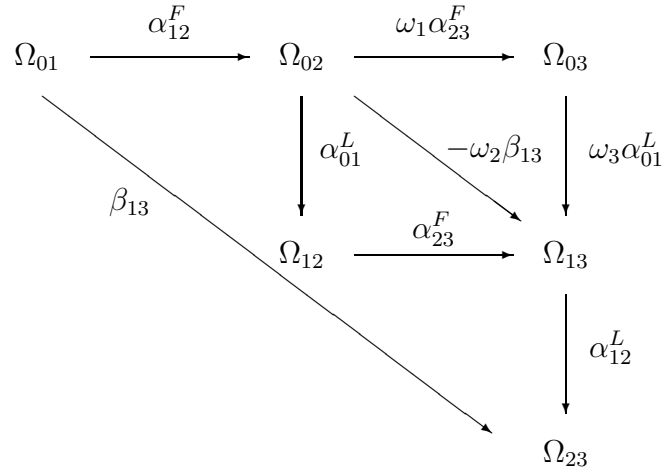
\begin{picture}(400,200)
\put(50,175){\makebox(0,0){$\om_{01} $}}

\put(52,160){\vector(4,-3){170}}

\put(85,120){\makebox(0,0){$\cextiii_{13}$}}

\put(100,188){\makebox(0,0){$\cextiif_{12}$}}
\put(70,175){\vector(1,0){60}}
\put(150,175){\makebox(0,0){$\om_{02} $}}
\put(170,175){\vector(1,0){60}}
\put(200,188){\makebox(0,0){$\k_1\cextiif_{23}$}}

\put(170,160){\vector(4,-3){63}}

\put(200,110){\makebox(0,0){$ \cextiif_{23}$}}

\put(150,160){\vector(0,-1){45}}
\put(165,140){\makebox(0,0){$\cextiil_{01}$}}
\put(223,140){\makebox(0,0){$-\k_2 \cextiii_{13}$}}

\put(273,140){\makebox(0,0){$\k_3\cextiil_{01}$}}
\put(267,65){\makebox(0,0){$\cextiil_{12}$}}

\put(250,160){\vector(0,-1){45}}
\put(250,85){\vector(0,-1){45}}

\put(250,175){\makebox(0,0){$\om_{03} $}}
\put(150,100){\makebox(0,0){$\om_{12} $}}

\put(170,100){\vector(1,0){60}}

\put(250,100){\makebox(0,0){$\om_{13} $}}
\put(250,25){\makebox(0,0){$\om_{23} $}}
\end{picture}

We recall that up to now we have referred in this paper to the cohomology
groups of Lie {\it algebras} and not of groups.
As mentioned they do not necessarily coincide \cite{Bar,VanEst,TUYNMAN},
as illustrated by the standard example of the $G$(2+1) Galilei group
(our $SO_{(0,0,1)}(4)$)
for which dim[$H^2(G(2+1),U(1))$]=2 although
dim[$H^2(\G(2+1),\Re)$]=3, a fact known
\cite{LL} (see also \cite{DoeMan} for a recent discussion)
for already twenty-five years\footnote{
This well known result has attracted a renewed interest \cite{Grigore},
specially in relation with the absence of non-relativistic planar systems with
exotic angular momentum (anyons) \cite{Bose}.}.
With $\om_{01} \equiv H$ and space rotation generator $\om_{23} \equiv J$,
we see that the algebra commutator
$[H, J]$ admits an extension through
$\cextiii_{13} $
but that the compactness condition on the space rotation
generator
$J$, relevant for the Galilei group, forces the coefficient
$\cextiii_{13}$ to disappear.
This is because under a rotation generated by $J$ in the
extended algebra, $H$ transforms by
$H \to \exp{(\theta J)}H\exp{(-\theta J)} = H - \theta
\cextiii_{13} \Xi$.
Since $J$ is compact the rotations
$\theta=2\pi$ (with $\k_3=1$) and
$\theta=0$ should coincide, which forces
$\cextiii_{13}=0$ and reduces in one dimension the group cohomology.
In general, within the CK
family of groups as obtained by exponentiation of the matrix
representation (\ref{ccd}) of the  CK Lie algebra, the one-parameter
subgroup generated by
$\om_{ab}$ is compact if $\k_{ab} > 0$ and non compact otherwise.
This implies that the extension
coefficient $\cextiii_{a+1\,c+1}$, which appears in the extended commutator
$[\om_{a\,a+1}, \om_{c\,c+1}] =\cextiii_{a+1\,c+1} \Xi$ for the algebra,
does not
correspond to a  group extension whenever at least one of the generators
$\om_{a\,a+1}$ or $\om_{c\,c+1}$ corresponds to
a compact one-parameter subgroup,
that is when either $\k_{a+1} > 0$ or $\k_{c+1} > 0$.


\subsect{$\overline{so}_{\k_1,\k_2,\k_3,\k_4}(5)$}

There are six basic extension coefficients of type II and three of type
III:
\bea
 && \cextiil_{01}\equiv\cext_{02,12} \qquad
 \cextiil_{12}\equiv \cext_{13,23}\qquad
 \cextiil_{23}\equiv \cext_{24,34}\cr
 && \cextiif_{12}\equiv \cext_{01,02}
\qquad  \cextiif_{23}\equiv\cext_{12,13} \qquad
\cextiif_{34}\equiv \cext_{23,24} \cr
 &&\cextiii_{13}\equiv \cext_{01,23}\qquad
\cextiii_{14}\equiv\cext_{01,34}\qquad
\cextiii_{24}\equiv\cext_{12,34}
 \label{fj}
\eea
verifying the additional conditions:
\be
\begin{array}{lll}
\k_3\cextiif_{12}  =\k_1\cextiil_{12}    &
 \quad  \k_4\cextiif_{23}
=\k_2\cextiil_{23}   &\quad    \cr
\k_1\k_2\cextiii_{13}=0  & \quad
\k_2\k_3\cextiii_{13}=0   &\quad
 \k_4\cextiii_{13}=0   \cr
 \k_2\cextiii_{14}=0   & \quad
 \k_3\cextiii_{14}=0  &\quad  \cr
\k_1\cextiii_{24}=0   &  \quad
\k_2\k_3\cextiii_{24}=0   & \quad
\k_3\k_4\cextiii_{24}=0\quad.
\end{array}
\label{fk}
\ee
Therefore the Lie brackets of the extended CK algebras
$\overline{so}_{\k_1,\k_2,\k_3,\k_4}(5)$ are:
\bea
&&[\om_{01},\om_{02}] =\k_1\om_{12}+\cextiif_{12}\Xi\qquad\qquad\qquad
[\om_{02},\om_{12}] =\k_2\om_{01}+\cextiil_{01}\Xi\cr
&&[\om_{12},\om_{13}] =\k_2\om_{23}+\cextiif_{23}\Xi\qquad\qquad\qquad
[\om_{13},\om_{23}] =\k_3\om_{12}+\cextiil_{12}\Xi\cr
&&[\om_{23},\om_{24}] =\k_3\om_{34}+\cextiif_{34}\Xi\qquad\qquad\qquad
[\om_{24},\om_{34}] =\k_4\om_{23}+\cextiil_{23}\Xi\cr
&&[\om_{13},\om_{14}] =\k_2(\k_3\om_{34}+\cextiif_{34}\Xi)\qquad\qquad
[\om_{04},\om_{14}] =\k_3\k_4(\k_2\om_{01}+\cextiil_{01}\Xi)\cr
 &&[\om_{02},\om_{03}] =\k_1(\k_2\om_{23}+\cextiif_{23}\Xi)\qquad\qquad
[\om_{14},\om_{24}] =\k_4(\k_3\om_{12}+\cextiil_{12}\Xi)\cr
&&[\om_{03},\om_{04}] =\k_1\k_2(\k_3\om_{34}+\cextiif_{34}\Xi)\qquad\quad
[\om_{03},\om_{13}] =\k_3(\k_2\om_{01}+\cextiil_{01}\Xi)\cr
&&[\om_{01},\om_{23}] =\cextiii_{13}\Xi\qquad
[\om_{02},\om_{13}] =-\k_2\cextiii_{13}\Xi\qquad
[\om_{01},\om_{34}] =\cextiii_{14}\Xi\cr
&&[\om_{12},\om_{34}] =\cextiii_{24}\Xi
 \qquad [\om_{13},\om_{24}] =-\k_3\cextiii_{24}\Xi\qquad,
\label{fm}
\eea
the remaining commutators being as in the non-extended case
(\ref{cb}).

We display the explicit result for  each CK algebra
${so}_{\k_1,\k_2,\k_3,\k_4}(5)$ in Table II. The first column shows the
number of simple contractions  (the number of coefficients $\k_a$ set equal
to zero).
The second schematically names the centrally extended Lie algebras.
The third specifies the coefficients $\k_a$ different
from zero together with the  non-trivial central extension coefficients
allowed.  Finally, the fourth gives
dim[$H^2({so}_{\k_1,\k_2,\k_3,\k_4}(5),\Re)$]
as a sum of the type II and type III contributions.  Note
that the only kinematical algebras in $(3+1)$ dimensions which  have
non-trivial central extensions (and hence projective representations)
are the (3+1) oscillating  Newton--Hooke  $(1,0,1,1)$, expanding
Newton--Hooke $(-1,0,1,1)$ and  Galilean $(0,0,1,1)$ algebras,  all of
them of  `absolute time' \cite{BLL}.  This table can be used as an
example of how
to compute dim[$H^2(\G,\Re)$] from Th. \teor\ .

\goodbreak

\bigskip
\noindent
\begin{tabular}{|l|l|l|l|}
\hline
\multicolumn{1}{|c|}{$\!$\#$\!$}&\multicolumn{1}{c|}{\em Extended algebra}&
\multicolumn{1}{c|}{\em (CK constants) [Non-trivial ext. coefficients]$\!$}&
\multicolumn{1}{c|}{\em $\!$dim$H^2\!$}\\
\hline
0&$\overline{ {so}}(5)$&($\k_1,\k_2,\k_3,\k_4$)&  0\\
 & $\overline{ {so}} (4,1)$& & \cr
 &$\overline{ {so}}(3,2)$ & & \\
\hline
&$\overline{ {iso}}(4)$&($0,\k_2,\k_3,\k_4$)
or ($\k_1,\k_2,\k_3,0$)& 0\\
 &$\overline{ {iso}}(3,1)$ & & \cr
 &$\overline{ {iso}}(2,2)$ & & \cr
1& & &\cr
&$\overline{ {t}}_6( {so}(3)\oplus
 {so}(2))$ &($\k_1,0,\k_3,\k_4$) \qquad [$\cextiil_{01}$]
or &1+0\cr
&$\overline{ {t}}_6( {so}(3)\oplus   {so}(1,1))$
 & ($\k_1,\k_2,0,\k_4$) \qquad [$\cextiif_{34}$]&\cr
&$\overline{ {t}}_6( {so}(2,1)\oplus   {so}(2))$ & &\cr
&$\overline{ {t}}_6( {so}(2,1)\oplus   {so}(1,1))$ & &\cr
\hline
&$\overline{ {iiso}}(3)$ &($0,0,\k_3,\k_4$) \qquad [$\cextiil_{01}$]
or &1+0\cr
&$\overline{ {iiso}}(2,1)$ &($\k_1,\k_2,0,0$) \qquad [$\cextiif_{34}$] &\cr
& & &\cr
&$\overline{ {ii'so}}(3)$ &($0,\k_2,\k_3,0$) &0\cr
&$\overline{ {ii'so}}(2,1)$ & &\cr
2& & &\cr
&$\overline{ {it}}_4( {so}(2)\oplus   {so}(2))$
&
($0,\k_2,0,\k_4$) \qquad [$\cextiil_{12},\cextiif_{12},
\cextiif_{34};\cextiii_{24}$] or &3+1\cr
&$\overline{ {it}}_4( {so}(2)\oplus   {so}(1,1))$
 &($\k_1,0,\k_3,0$) \qquad [$\cextiil_{01},\cextiil_{23},
\cextiif_{23};\cextiii_{13}$] &\cr
&$\overline{ {it}}_4( {so}(1,1)\oplus   {so}(1,1))$ & &\cr
& & &\cr
& $\overline{ {t}}_6( {iso}(2)\oplus  {so}(2))$&
($\k_1,0,0,\k_4$) \qquad [$\cextiil_{01},\cextiif_{34};\cextiii_{14}$]
 &2+1\cr
& $\overline{ {t}}_6( {iso}(2)\oplus  {so}(1,1))$& &\cr
& $\overline{ {t}}_6( {iso}(1,1)\oplus   {so}(1,1))$& &\cr
\hline
&$\overline{ {iiiso}}(2)$ &
($0,0,0,\k_4$) \qquad [$\cextiil_{01},
\cextiif_{12},\cextiil_{12},\cextiif_{34};\cextiii_{14},\cextiii_{24}$]
or &4+2\cr
3&$\overline{ {iiiso}}(1,1)$ &
($\k_1,0,0,0$) \qquad [$\cextiil_{01},
\cextiif_{23},\cextiil_{23},
\cextiif_{34},\cextiii_{13},\cextiii_{14}$]&\cr
& & &\cr
& $\overline{ {iii'so}}(2)$&
($0,0,\k_3,0$) \qquad [$\cextiil_{01},
\cextiif_{23},\cextiil_{23};\cextiii_{13},\cextiii_{24}$]
or&3+2\cr &$\overline{
{iii'so}}(1,1)$ &
($0,\k_2,0,0$) \qquad [$\cextiif_{12},
\cextiil_{12},\cextiif_{34};\cextiii_{13},\cextiii_{24}$] &\cr
\hline
4& $\overline{ {iiiiso}}(1)$& ($0,0,0,0$)&6+3\cr
&&\qquad [$\cextiil_{01},\cextiif_{12},
\cextiil_{12},\cextiif_{23},\cextiil_{23},\cextiif_{34};
\cextiii_{13},\cextiii_{14},\cextiii_{24}$]
&\\
\hline
\end{tabular}

\smallskip
{\noindent\centerline{
{\bf {Table II}}. Non-trivial central extensions
$\overline{so}_{\k_1,\k_2,\k_3,\k_4}(5)$ of ${so}_{\k_1,\k_2,\k_3,\k_4}(5)$.}
\noindent\centerline{
The constants $\k_i$ appearing explicitly are assumed to be different from
zero.
}}

\baselineskip=18pt


\sect {Concluding remarks}

We have characterised with
generality  the
second cohomology groups \allowbreak
$H^2(so_{\k_1\ldots\k_N}(N+1),\Re)$
of the CK family of algebras
$so_{\k_1\ldots\k_N}(N+1)$,
which is a particular subfamily of all graded
contractions of the $so(N+1)$ algebra.
The algebras in the CK family can be described in a
simultaneous and economical way using $N$ real `contraction' coefficients
$\k_1, \k_2, \dots, \k_N$.
The procedure also exhibits the origin of the
various central extensions, and in particular differentiates clearly those
which come from contractions of trivial extensions from those which do not.

It is well known that, by Whitehead's lemma, all semisimple Lie
algebras have trivial  second
cohomology groups, and that by the
Levi--Mal'\v{c}ev theorem any finite  dimensional Lie algebra $\G$ is
the semidirect extension of a semisimple  algebra and the radical of
$\G$. Since inhomogeneous algebras come from contraction, our
procedure may be applied to find the cohomology groups of
other inhomogeneous algebras as well; in particular, one could start from
the real simple algebras of the $A_l$ and $C_l$ series.
There are several CK families of
algebras (see \cite{HerSanGoslar96} for a cursory description), and {\em
any} simple real Lie algebra appears as a member of some family. We have
discussed here only the orthogonal CK family, which include the simple
algebras $so(N+1)$ and $so(p,q)$ in the $B_l$ and $D_l$ series, as well as
their
(quasi-simple) contractions.
A similar approach would lead to a
complete characterisation of the second cohomology groups for
quasi-simple algebras of inhomogeneous type obtained  by contraction from
other real simple Lie algebras.
This may be the matter for further work.

Another possible application of the contraction method
is the search for  Casimir operators of inhomogeneous
algebras. The  number of primitive Casimirs of a simple algebra
$\G$ is  equal to its rank
$l$, which in turn is equal to the different  primitive invariant
polynomials which can be constructed on $\G$.
Thus, the graded
contraction approach allows us, in principle, to find central  elements
of the enveloping algebras by contracting the original $l$
Casimir-Racah operators.
Clearly, the procedure does not permit to find {\it all} the Casimirs
of an {\it arbitrary} contraction of a simple Lie algebra of rank $l$,
since the final step is always an Abelian algebra (hence with as many
primitive Casimirs as generators) and dim$\G >l$.
However, within the CK family the number of functionally independent
Casimirs remains constant (see \cite{HS}).
This provides another justification for the name `quasi-simple'
given to its members, and explains in a simple way why {\it e.g.}, the number
of Casimir operators for the simple de Sitter algebra and the non-simple
Poincar\'e one is the same.

The same kind of approach we have pursued here for studying the
second cohomology groups of the
CK algebras has been developed
to study their deformations (in the sense of \cite{DRI,JIM,FRT}).
In particular, a whole family of deformations of
inhomogeneous Lie algebras \cite{BHOSd}, or working to first order, of the
corresponding bialgebras \cite{BGHOS}, has been found.
The semidirect structure of the `classical' CK $\k_1 = 0$ inhomogeneous Lie
algebras becomes \cite{bicross} a
bicrossproduct \cite{Majid} structure for their CK deformed counterparts.
Whether or not this extends to the deformations of other semidirect
structures associated to the vanishing of any $\k_a$ requires
further study.
A related problem would be the analysis of
the structure of the deformation of inhomogeneous Lie algebras from the present
graded contraction point of view, for which central extensions should appear
as cocycle-bicrossproducts. These questions are worth studying.

\section*{Acknowledgements}
This research has been partially supported by the Spanish DGES
projects PB96--0756, PB94--1115 and
PR95--439.  J.A. and J.C.P.B. wish to thank the kind hospitality extended to
them at  DAMTP. Finally, the support of St. John's College (J.A.) and
an FPI grant from the Spanish Ministry of Education and Science and the CSIC
(J.C.P.B.) are gratefully acknowledged.


\end{document}